\documentclass[12pt]{article}
\usepackage{latexsym, amssymb, amsmath, amscd, amsfonts, epsfig, graphicx, colordvi,verbatim,ifpdf,extarrows}
\usepackage{amsfonts, amsmath, amssymb}
\usepackage{amssymb,amsfonts,amsmath,latexsym,epsfig,cite, psfrag,eepic,color}
\usepackage{amscd,graphics}
\usepackage{latexsym, amssymb,  amsmath,amscd, amsfonts, epsfig, graphicx, colordvi,amsthm}

\usepackage{graphicx}
\usepackage{epstopdf}
\usepackage{color}
\usepackage{ifpdf}
\usepackage{fancybox}
\usepackage[font=small,labelfont=bf,labelsep=none]{caption}
\usepackage{float}

\newtheorem{thm}{Theorem}[section]
\newtheorem{pro}[thm]{Proposition}

\newtheorem{defi}[thm]{Definition}
\newtheorem{lem}[thm]{Lemma}

\def\pf{\noindent{\it Proof.} }
\def\qed{\nopagebreak\hfill{\rule{4pt}{7pt}}
\medbreak}

\numberwithin{equation}{section}

\def\qed{\nopagebreak\hfill{\rule{4pt}{7pt}}
\medbreak}

\def \GG {G\"ollnitz-Gordon marking}

\newlength{\boxedparwidth}
  {\begin{center} \begin{tabular}{|@{\hspace{.315in}}c@{\hspace{.15in}}|}
                  \hline \\ \begin{minipage}[t]{\boxedparwidth}
                  \setlength{\parindent}{.25in}}%
  {\end{minipage} \\ \\ \hline \end{tabular} \end{center}}

\begin{document}

\begin{center}

 {\Large \bf An overpartition analogue of \\[5pt] the Andrews-G\"ollnitz-Gordon theorem}
\end{center}

\begin{center}
 {Thomas Y. He}$^{1}$, {Kathy Q. Ji}$^{2}$, {Allison Y.F. Wang}$^{3}$ and
  {Alice X.H. Zhao}$^{4}$ \vskip 2mm

$^{1,2,3}$Center for Applied Mathematics,  Tianjin University, Tianjin 300072, P.R. China\\[6pt]

  $^{4}$ College of Science, Tianjin University of Technology, Tianjin 300384, P.R. China

   \vskip 2mm

    $^1$heyao@tju.edu.cn, $^2$kathyji@tju.edu.cn,  $^3$wangyifang@tju.edu.cn, $^4$zhaoxh@email.tjut.edu.cn
\end{center}

\vskip 6mm \noindent {\bf Abstract.} In 1967, Andrews found a combinatorial generalization
of the  G\"ollnitz-Gordon theorem, which can be called the Andrews-G\"ollnitz-Gordon theorem. In 1980, Bressoud derived a multisum Rogers-Ramanujan-type identity, which can be considered as the generating function counterpart of the Andrews-G\"ollnitz-Gordon theorem. Lovejoy gave an overpartition analogue of the  Andrews-G\"ollnitz-Gordon theorem for   $i=k$. In this paper, we give an overpartition analogue of this theorem in the general case.   By using  Bailey's lemma and a change of base formula due to Bressoud, Ismail and Stanton, we  obtain an overpartition analogue of Bressoud's identity. We then give a combinatorial interpretation of this identity  by introducing the G\"ollnitz-Gordon marking of an overpartition, which yields an overpartition analogue of the Andrews-G\"ollnitz-Gordon theorem.



\section{Introduction}

The purpose of this paper is to  give an overpartition analogue of  the Andrews-G\"ollnitz-Gordon theorem in the general case. In 1967, Andrews \cite{Andrews-1967} found the following  combinatorial generalization of the G\"ollnitz-Gordon identities \cite{Gollnitz-1967, Gordon-1965}, which has been called the   Andrews-G\"ollnitz-Gordon theorem.

\begin{thm}[Andrews-G\"ollnitz-Gordon]\label{Gordon-Gollnitz-odd}
For $k\geq i\geq1$, let $C_{k,i}(n)$ denote the number of partitions $\lambda$ of $n$ of the form $(1^{f_1},\,2^{f_{2}},\,3^{f_3},\ldots)$,  where $f_t(\lambda)$ {\rm(}or $f_t$ for short{\rm )} denotes the number of occurrences of  $t$ in   $\lambda$,  such that
\begin{itemize}
\item[{\rm (1)}] $f_1(\lambda)+f_2(\lambda)\leq i-1${\rm;}

\item[{\rm (2)}] $f_{2t+1}(\lambda)\leq1${\rm;}

\item[{\rm (3)}] $f_{2t}(\lambda)+f_{2t+1}(\lambda)+f_{2t+2}(\lambda)\leq k-1$.
\end{itemize}

For $k\geq i\geq1$, let $D_{k,i}(n)$ denote the number of partitions of $n$ into parts $\not\equiv2\pmod4$ and $\not\equiv0,\pm(2i-1)\pmod{4k}$.

Then for $k\geq i\geq1$ and  $n\geq0$,
\[C_{k,i}(n)=D_{k,i}(n).\]
\end{thm}

The Andrews-G\"ollnitz-Gordon theorem was motivated by a combinatorial generalization of the Rogers-Ramanujan identities due to Gordon \cite{Gordon-1961}.

\begin{thm}[Rogers-Ramanujan-Gordon]\label{R-R-G}
For $k\geq i\geq1$, let $B_{k,i}(n)$ denote the number of partitions $\lambda$ of $n$ of the form $(1^{f_1},\,2^{f_{2}},\,3^{f_3},\ldots)$ such that
\begin{itemize}
\item[{\rm(1)}] $f_1(\lambda)\leq i-1${\rm;}
\item[{\rm (2)}] $f_t(\lambda)+f_{t+1}(\lambda)\leq k-1$.
\end{itemize}
For $k\geq i\geq1$, let $A_{k,i}(n)$ denote the number of partitions of $n$ into parts $\not\equiv0,\pm i\pmod{2k+1}$.

Then for  $k\geq i\geq1$ and  $n\geq0$,
\[A_{k,i}(n)=B_{k,i}(n).\]
\end{thm}

The analytic proof of Theorem \ref{R-R-G} was provided by Andrews  \cite{Andrews-1966}, and in  \cite{Andrews-1974}, he discovered the following identity, which has been called the Andrews-Gordon identity, see Kur\c{s}ung\"{o}z \cite{Kursungoz-2010}.

\begin{thm}[Andrews]\label{R-R-G-e}
For $k\geq i\geq 1$,
\begin{equation}\label{R-R-A1}
\sum_{N_1\geq N_2\geq \ldots \geq N_{k-1}\geq0}\frac{q^{N_1^2+N_2^2+\cdots+N^2_{k-1}+N_i+\cdots+N_{k-1}}}
{(q;q)_{N_1-N_2}(q;q)_{N_2-N_3}\cdots(q;q)_{N_{k-1}}}
=\frac{(q^i,q^{2k+1-i},q^{2k+1};q^{2k+1})_\infty}{(q;q)_\infty}.
\end{equation}
\end{thm}
Here and in the sequel, we  adopt the standard notation \cite{Andrews-1976}:
\[(a;q)_\infty=\prod_{i=0}^{\infty}(1-aq^i), \quad (a;q)_n=\frac{(a;q)_\infty}{(aq^n;q)_\infty},\]
and
\[(a_1,a_2,\ldots,a_m;q)_\infty=(a_1;q)_\infty(a_2;q)_\infty\cdots(a_m;q)_\infty.
\]

 Theorem \ref{R-R-G-e}    can be considered as the generating function version of Theorem \ref{R-R-G}. It is evident that the generating function of $A_{k,i}(n)$  in Theorem \ref{R-R-G}  equals the right-hand side of \eqref{R-R-A1}. By using   $q$-difference equations, Andrews  \cite{Andrews-1974} showed that the generating function of $B_{k,i}(n)$   in Theorem \ref{R-R-G} equals the left-hand side of \eqref{R-R-A1}. In particular, he obtained the following formula for the generating function of   $B_{k,i}(m,n)$, where  $B_{k,i}(m,n)$ denotes the number of partitions enumerated by $B_{k,i}(n)$ with exactly $m$ parts.
\begin{thm} [Andrews]\label{R-R-G-left-e}
For $k\geq i\geq 1$,
\begin{equation}\label{R-R-A1-gf}
\sum_{m,\,n\geq 0}B_{k,i}(m,n)x^mq^n=\sum_{N_1\geq N_2\geq \ldots \geq N_{k-1}\geq0}\frac{q^{N_1^2+N_2^2+\cdots+N^2_{k-1}+N_i
+\cdots+N_{k-1}}x^{N_1+\cdots+N_{k-1}}}{(q;q)_{N_1-N_2}
(q;q)_{N_2-N_3}\cdots(q;q)_{N_{k-1}}}.
\end{equation}
\end{thm}

Kur\c{s}ung\"{o}z \cite{Kursungoz-2010} gave a combinatorial proof of \eqref{R-R-A1-gf} by introducing the notion of the Gordon marking of a partition.

The generating function version of the Andrews-G\"ollnitz-Gordon theorem \ref{Gordon-Gollnitz-odd} was found by Bressoud \cite[Eq. (3.8)]{Bressoud-1980}  in 1980.

\begin{thm}[Bressoud]\label{Gordon-Gollnitz-odd-e}
For $k\geq i\geq1$,
\begin{equation}\label{Gollnitz-A1}
\begin{split}
&\sum_{N_1\geq N_2\geq \ldots \geq N_{k-1}\geq0}\frac{(-q^{1-2N_1};q^2)_{N_1}
q^{2(N^2_1+\cdots+N^2_{k-1}+N_i+\cdots+N_{k-1})}}
{(q^2;q^2)_{N_1-N_2}(q^2;q^2)_{N_2-N_3}\cdots(q^{2};q^{2})_{N_{k-1}}}\\
&=\frac{(q^2;q^4)_\infty(q^{4k},q^{2i-1},q^{4k-2i+1};q^{4k})_\infty}{(q;q)_\infty}.
\end{split}
\end{equation}
\end{thm}

 Bressoud \cite{Bressoud-1980} also showed that  the left-hand side of \eqref{Gollnitz-A1} can be interpreted combinatorially as the generating function of $C_{k,i}(n)$  in Theorem \ref{Gordon-Gollnitz-odd}. More precisely, he gave the following formula for the generating function of $C_{k,i}(m,n)$, where $C_{k,i}(m,n)$ denotes the number of partitions enumerated by $C_{k,i}(n)$ with exactly $m$ parts.

\begin{thm}[Bressoud]\label{GENERATING-2}
For $k\geq i\geq1$,
\begin{equation}\label{GENERATING2}
\begin{split}
&\displaystyle\sum_{m,n\geq0}C_{k,i}(m,n)x^mq^n\\
&=\sum_{N_1\geq N_2\geq \ldots \geq N_{k-1}\geq0}\frac{(-q^{1-2N_1};q^2)_{N_1}
q^{2(N^2_1+\cdots+N^2_{k-1}+N_i+\cdots+N_{k-1})}x^{N_1+\cdots+N_{k-1}}}
{(q^2;q^2)_{N_1-N_2}(q^2;q^2)_{N_2-N_3}\cdots(q^{2};q^{2})_{N_{k-1}}}.
\end{split}
\end{equation}
\end{thm}

 In recent years,   many overpartition analogues of classical partition theorems  have been proved, see, for example, Chen, Sang and Shi \cite{Chen-Sang-Shi-2013}, Corteel and Mallet \cite{Corteel-2007}, Corteel, Lovejoy and Mallet \cite{Corteel-2008}, and Lovejoy \cite{Lovejoy-2003, Lovejoy-2004, Lovejoy-2007,Lovejoy-2010}. In particular, Lovejoy \cite{Lovejoy-2004} obtained an overpartition analogue of the  Andrews-G\"ollnitz-Gordon theorem for  $i=k$. In this paper, we give an overpartition analogue of the Andrews-G\"ollnitz-Gordon theorem  in the general case. We also obtain an overpartition analogue of  Bressoud's identity \eqref{Gollnitz-A1}.

Recall that an overpartition  of  $n$ is a partition of $n$ in which the first occurrence of a number can be overlined.  For an overpartition $\lambda$ of $n$, let $f_t(\lambda)$ (resp. $f_{\bar{t}}(\lambda)$) be the number of occurrences of    $t$ (resp. $\bar{t}$) in $\lambda$, without  ambiguity, we  write  $f_t$ (resp. $f_{\bar{t}}$) for short. By the definition of an overpartition, it is clear to see that $f_{\bar{t}}(\lambda)=0$ or $1$.

We   obtain the following overpartition analogue of the Andrews-G\"ollnitz-Gordon theorem.

\begin{thm}\label{Gordon-Gollnitz-odd-overp}
For $k\geq i\geq1$, let $O_{k,i}(n)$ denote the number of overpartitions of $n$ of the form $(\bar{1}^{f_{\bar{1}}},\,1^{f_1},\,\bar{2}^{f_{\bar{2}}},\,2^{f_2},\ldots)$  such that
\begin{itemize}
\item[{\rm(1)}] $f_{\overline{1}}(\lambda)+f_2(\lambda)\leq i-1${\rm;}

\item[{\rm (2)}] $f_{\overline{2t}}(\lambda)+ f_{2t}(\lambda)+f_{\overline{2t+1}}(\lambda)+f_{2t+2}(\lambda)\leq k-1${\rm;}
\item[{\rm (3)}]If $f_{2t+1}(\lambda)\geq1$, then $f_{2t+2}(\lambda)\leq k-2$.
\end{itemize}
For $k>i\geq 1$, let $P_{k,i}(n)$ denote the number of overpartitions of $n$ with non-overlined parts   $\not\equiv 0,\pm (2i-1)$  $\pmod{4k-2}$ and let $P_{k,k}(n)$ denote the number of overpartitions of $n$ into parts not divisible by $2k-1$.

Then for  $k\geq i\geq1$ and  $n\geq0$,
\[O_{k,i}(n)=P_{k,i}(n).\]
 \end{thm}

 It should be noted that for an overpartition $\lambda$ counted by $O_{k,i}(n)$ without overlined even parts and non-overlined odd parts,  if we change  overlined odd parts of $\lambda$ to non-overlined odd parts, then we get a partition  counted by $C_{k,i}(n)$. Hence we could say that Theorem \ref{Gordon-Gollnitz-odd-overp} is an overpartition analogue of Theorem \ref{Gordon-Gollnitz-odd}.

 We also obtained the following overpartition analogue of Bressoud's identity \eqref{Gollnitz-A1}, which can be viewed as the corresponding generating function version of Theorem \ref{Gordon-Gollnitz-odd-overp}.

\begin{thm}\label{Gordon-Gollnitz-odd-overp-eqn}
For $k\geq i\geq1$,
\begin{align}\label{over3}
&\sum_{N_{1}\geq \cdots\geq N_{k-1}\geq 0}\frac{(-q^{2-2N_1};q^2)_{N_1-1}(-q^{1-2N_1};q^2)_{N_1}
q^{2(N^2_1+\cdots+N^2_{k-1}+N_{i+1}+\cdots+N_{k-1})}
(1+q^{2N_i})}{(q^2;q^2)_{N_1-N_2}\cdots(q^2;q^2)_{N_{k-2}-N_{k-1}}
(q^{2};q^{2})_{N_{k-1}}}\nonumber\\
&=\frac{(-q;q)_\infty(q^{2i-1},q^{4k-1-2i},q^{4k-2};q^{4k-2})_\infty}
{(q;q)_\infty}.
\end{align}
\end{thm}

We will first  prove Theorem \ref{Gordon-Gollnitz-odd-overp-eqn}  by using Bailey's lemma and a  change of base formula due to   Bressoud, Ismail and Stanton \cite{Bressoud-Ismail-Stanton-2000}. We then use  Theorem \ref{Gordon-Gollnitz-odd-overp-eqn} to derive Theorem \ref{Gordon-Gollnitz-odd-overp}. More precisely,  let $O_{k,i}(m,n)$ denote the number of overpartitions counted by $O_{k,i}(n)$ with exactly $m$ parts, we shall give a combinatorial proof of the following formula for the  generating function  of $O_{k,i}(m,n)$ by introducing the G\"ollnitz-Gordon marking of an overpartition.

\begin{thm}\label{Gollnitz-odd-1}
For $k\geq i\geq1$,
\begin{equation}\label{Gollnitz-odd1}
\begin{split}
& \displaystyle\sum_{m,n\geq0}O_{k,i}(m,n)x^mq^n \\[5pt]
&=\sum_{N_{1}\geq \cdots\geq N_{k-1}\geq 0}\frac{(-q^{2-2N_1};q^2)_{N_1-1}(-q^{1-2N_1};q^2)_{N_1}}
{(q^2;q^2)_{N_1-N_2}\cdots}\\[5pt]
&\quad\quad \quad \quad \quad  \times \frac{q^{2(N^2_1+\cdots+N^2_{k-1}+N_{i+1}+\cdots+N_{k-1})}(1+q^{2N_i})
x^{N_1+\cdots+N_{k-1}}}{(q^2;q^2)_{N_{k-2}-N_{k-1}}(q^{2};q^{2})_{N_{k-1}}}.
\end{split}
\end{equation}
\end{thm}

Setting $x=1$ in \eqref{Gollnitz-odd1}, we obtain the generating function for $O_{k,i}(n)$ which is the left-hand side of  \eqref{over3}. On the other hand, it is evident that the generating function of $P_{k,i}(n)$ equals
\begin{equation}\label{Gollnitz-odd-left}
\displaystyle\sum_{n\geq0}P_{k,i}(n)q^n=\frac{(-q;q)_\infty(q^{2i-1},q^{4k-2i-1},q^{4k-2};q^{4k-2})_\infty}
{(q;q)_\infty},
\end{equation}
which is the right-hand side of  \eqref{over3}. Hence we are led to  Theorem \ref{Gordon-Gollnitz-odd-overp} by Theorem  \ref{Gordon-Gollnitz-odd-overp-eqn}.

 The paper is organized as follows. In Section 2, we first review some necessary results on Bailey pairs and then give a proof of Theorem  \ref{Gordon-Gollnitz-odd-overp-eqn} by combining Bailey's lemma and a change of base formula.  In Section 3, we begin with introducing the notion of the G\"ollnitz-Gordon marking of an overpartition  and  then give an outline of the proof of   the formula for the generating function
 of  $O_{k,i}(m,n)$  in Theorem \ref{Gollnitz-odd-1}. It turns out that  the proof of Theorem \ref{Gollnitz-odd-1} reduces to  the proofs of two relations stated in Lemma \ref{N-1} and Lemma  \ref{N-2}, respectively. Section 4 and Section 5 are devoted to the bijective proofs of these two relations respectively.  In Section 6,  we complete  the proof of Theorem \ref{Gollnitz-odd-1}.

\section{Proof of Theorem \ref{Gordon-Gollnitz-odd-overp-eqn}}

  We  first give a brief review  of some relevant results on  Bailey pairs required in the proof of Theorem \ref{Gordon-Gollnitz-odd-overp-eqn}.   For more on Bailey pairs, see, for example, \cite{Agarwal-Andrews-Bressoud, Andrews-1986, Andrews-2000, Bressoud-Ismail-Stanton-2000, Lovejoy-2004b, Paule-1987,Warnaar-2001}.  Recall that  a pair of sequences $(\alpha_n(a,q),\beta_n(a,q))$ is called a Bailey pair with parameters $(a,q)$ (or  a Bailey pair for short) if  for all $n\geq 0,$
\begin{equation}\label{bailey pair}
\beta_n(a,q)=\sum_{r=0}^n\frac{\alpha_r(a,q)}{(q;q)_{n-r}(aq;q)_{n+r}}.
\end{equation}
  Bailey's lemma was first given by Bailey \cite{Bailey-1949}  and was formulated by Andrews  \cite{Andrews-1984, Andrews-1986} in the following form.

\begin{thm}[Bailey's lemma]\label{bailey lemma}
If $(\alpha_n(a,q),\beta_n(a,q))$ is a Bailey pair,
then $(\alpha'_n(a,q),\beta'_n(a,q))$ is also a Bailey pair, where
\begin{equation}\label{bl}
\begin{split}
\alpha_n'(a,q)&=\frac{(\rho_1;q)_n(\rho_2;q)_n}{(aq/\rho_1;q)_n(aq/\rho_2;q)_n}\left(\frac{aq}
{\rho_1\rho_2}\right)^n\alpha_n(a,q),\\[5pt]
\beta_n'(a,q)&=\sum_{j=0}^{n}\frac{(\rho_1;q)_j(\rho_2;q)_j(aq/\rho_1\rho_2;q)_{n-j}}
{(aq/\rho_1;q)_n(aq/\rho_2;q)_n(q;q)_{n-j}}\left(\frac{aq}{\rho_1\rho_2}\right)^j\beta_j(a,q).
\end{split}
\end{equation}
\end{thm}

Andrews first noticed that Bailey's lemma can create a new Bailey pair from a given one.  Hence iterating Theorem \ref{bailey lemma} produces a sequence of Bailey pairs, which has been called a Bailey chain.  Based on this observation,  Andrews \cite{Andrews-1984} showed that the Andrews-Gordon identity \eqref{R-R-A1} in Theorem \ref{R-R-G-e}  holds for $i=1$ and $i=k$ by iteratively using the following  specialization of  Bailey's lemma.

\begin{lem}[$\rho_1,\rho_2\rightarrow \infty$ in Theorem \ref{bailey lemma}]  \label{bl1}
If $(\alpha_n(a,q),\beta_n(a,q))$ is a Bailey pair,
then $(\alpha_n'(a,q), \beta_n'(a,q))$ is also a Bailey pair, where
\begin{equation}
\begin{split}
&\alpha_n'(a,q)=a^nq^{n^2}\alpha_n(a,q),\\[3pt]
&\beta_n'(a,q)=\sum_{j=0}^n\frac{a^jq^{j^2}}{(q;q)_{n-j}}\beta_j(a,q).
\end{split}
\end{equation}
\end{lem}

Subsequently, Agarwal, Andrews and Bressoud \cite{Agarwal-Andrews-Bressoud} gave an extension of a Bailey chain known as a Bailey lattice and used the Bailey lattice to  prove      the Andrews-Gordon identity \eqref{R-R-A1}  holds for $1\leq i\leq k$.  Bressoud, Ismail and Stanton \cite{Bressoud-Ismail-Stanton-2000} provided an alternative proof of  this identity    by  combining Bailey's lemma with the following proposition.

\begin{pro}\label{bl4}{\rm  \cite[Proposition 4.1]{Bressoud-Ismail-Stanton-2000}}
If $(\alpha_n(1,q),\beta_n(1,q))$ is a Bailey pair, where
\[\alpha_n(1,q)=\left\{
             \begin{array}{ll}
               1, & \hbox{for $n=0$}, \\[5pt]
               (-1)^nq^{An^2}(q^{(A-1)n}+q^{-(A-1)n}), & \hbox{for $n\geq 1$,}
             \end{array}
           \right.
\]
then $(\alpha_n'(1,q),\beta_n'(1,q))$ is also a Bailey pair, where
\begin{align*}
 \alpha_n'(1,q)&=\left\{
             \begin{array}{ll}
               1, & \hbox{for $n=0$}, \\[3pt]
               (-1)^nq^{An^2}(q^{An}+q^{-An}), & \hbox{for $n\geq 1$,}
             \end{array}
           \right. \\[10pt]
 \beta_n'(1,q)&=q^n\beta_n(1,q).
\end{align*}
\end{pro}

By iteratively using Bailey's lemma and Proposition \ref{bl4}, Bressoud, Ismail and Stanton \cite{Bressoud-Ismail-Stanton-2000} also provided a proof of Bressoud's identity \eqref{Gollnitz-A1} in Theorem \ref{Gordon-Gollnitz-odd-e}. Moreover, they   established new versions of Bailey's lemma, known as change of base formulas,  which change the base in Bailey pairs from $q$ to $q^2$ or $q^3$.  Iterating  these change of base formulas, they obtained many  new multisum Rogers-Ramanujan identities.

To prove Theorem \ref{Gordon-Gollnitz-odd-overp-eqn}, we shall make use of the following special case of the change of base formula due to Bressoud, Ismail and Stanton.

\begin{lem}\label{bailey lemma3} {\rm \cite[Theorem 2.5, $B\rightarrow \infty$]{Bressoud-Ismail-Stanton-2000}}
If $(\alpha_n(a,q),\beta_n(a,q))$ is a Bailey pair,
then $(\alpha_n'(a,q), \beta_n'(a,q))$ is also a Bailey pair, where
\begin{equation*}
\begin{split}
  \alpha_n'(a,q)&=\frac{1+a}{1+aq^{2n}}q^n\alpha_n(a^2,q^2),\\[3pt]
  \beta_n'(a,q)&=\sum_{k=0}^n\frac{(-a;q)_{2k}q^k}{(q^2;q^2)_{n-k}}\beta_{k}(a^2,q^2).
\end{split}
\end{equation*}
\end{lem}

We are now in a position to prove Theorem \ref{Gordon-Gollnitz-odd-overp-eqn}.

\noindent{\bf Proof of Theorem \ref{Gordon-Gollnitz-odd-overp-eqn}: }
We begin with the   unit Bailey pair \cite[H(17)]{Slater-1952},
\begin{equation}\label{BP1}
\begin{split}
&\alpha^{(0)}_n(1,q)=\left\{
             \begin{array}{ll}
               1, & \hbox{if $n=0$}, \\[5pt]
               (-1)^nq^{n^2/2}(q^{-n/2}+q^{n/2}), & \hbox{if $n\geq1$,}
             \end{array}
               \right.\\[5pt]
&\beta^{(0)}_n(1,q)=\left\{
             \begin{array}{ll}
               1, & \hbox{if $n=0$}, \\[5pt]
               0, & \hbox{if $n\geq1$.}
             \end{array}
               \right.
\end{split}
\end{equation}

Substituting  \eqref{BP1} into Lemma \ref{bl1}, we  get a new Bailey pair $(\alpha^{(1)}_n(1,q),\beta^{(1)}_n(1,q))$,
\begin{equation}\label{BP2}
\begin{split}
&\alpha^{(1)}_n(1,q)=\left\{
             \begin{array}{ll}
               1, & \hbox{if $n=0$}, \\[5pt]
               (-1)^nq^{3n^2/2}(q^{-n/2}+q^{n/2}), & \hbox{if $n\geq1$,}
             \end{array}
               \right.\\[5pt]
&\beta^{(1)}_n(1,q)=\frac{1}{(q;q)_n}.
\end{split}
\end{equation}

Alternatively applying Proposition \ref{bl4} and Lemma \ref{bl1}  $k-i-1$ times to \eqref{BP2} yields the following Bailey pair $(\alpha^{(2k-2i-1)}_n(1,q),\beta^{(2k-2i-1)}_n(1,q))$,
\begin{equation}\label{BP4}
\begin{split}
&\alpha^{(2k-2i-1)}_n(1,q)=\left\{
             \begin{array}{ll}
               1, & \hbox{if $n=0$}, \\[5pt]
               (-1)^nq^{\frac{2k-2i+1}{2}n^2}(q^{-\frac{2k-2i-1}{2}n}+q^{\frac{2k-2i-1}{2}n}), & \hbox{if $n\geq1$,}
             \end{array}
               \right.\\[10pt]
&\beta^{(2k-2i-1)}_n(1,q)
=\sum_{n\geq N_{i+1}\geq \cdots\geq N_{k-1}\geq 0}\frac{q^{N^2_{i+1}+N^2_{i+2}+\cdots+N^2_{k-1}
+N_{i+1}+\cdots+N_{k-1}}}{(q;q)_{n-N_{i+1}}(q;q)_{N_{i+1}-N_{i+2}}\cdots(q;q)_{N_{k-2}-N_{k-1}}(q;q)_{N_{k-1}}}.
\end{split}
\end{equation}
 Applying Proposition \ref{bl4}  to \eqref{BP4} gives the Bailey pair $(\alpha^{(2k-2i)}_n(1,q),\beta^{(2k-2i)}_n(1,q))$,
 \begin{equation}\label{BP5}
\begin{split}
&\alpha^{(2k-2i)}_n(1,q)=\left\{
             \begin{array}{ll}
               1, & \hbox{if $n=0$}, \\[5pt]
               (-1)^nq^{\frac{2k-2i+1}{2}n^2}(q^{-\frac{2k-2i+1}{2}n}+q^{\frac{2k-2i+1}{2}n}), & \hbox{if $n\geq1$,}
             \end{array}
               \right.\\[10pt]
&\beta^{(2k-2i)}_n(1,q)=q^n\sum_{n\geq N_{i+1}\geq \cdots\geq N_{k-1}\geq 0}\frac{q^{N^2_{i+1}+N^2_{i+2}+\cdots+N^2_{k-1}+N_{i+1}+\cdots+N_{k-1}}}{(q;q)_{n-N_{i+1}}(q;q)_{N_{i+1}-N_{i+2}}\cdots(q;q)_{N_{k-2}-N_{k-1}}(q;q)_{N_{k-1}}}.
\end{split}
\end{equation}
By the definition of a Bailey pair, it is easy to see that the sum of \eqref{BP4} and \eqref{BP5} generates a new Bailey pair $(\alpha^{(2k-2i+1)}_n(1,q),\,\beta^{(2k-2i+1)}_n(1,q))$,
\begin{equation}\label{BP6}
\begin{split}
&\alpha^{(2k-2i+1)}_n(1,q)=\left\{
             \begin{array}{ll}
               1, & \hbox{if $n=0$}, \\[5pt]
               (-1)^nq^{\frac{2k-2i+1}{2}n^2}(q^{-\frac{2k-2i+1}{2}n}
               +q^{\frac{2k-2i-1}{2}n})(1+q^n)/2, & \hbox{if $n\geq1$,}
             \end{array}
               \right.\\[10pt]
&\beta^{(2k-2i+1)}_n(1,q)=\sum_{n\geq N_{i+1}\geq \cdots\geq N_{k-1}\geq 0}\frac{(1+q^n) q^{N^2_{i+1}+N^2_{i+2}+\cdots+N^2_{k-1}+N_{i+1}+\cdots+N_{k-1}}}{2(q;q)_{n-N_{i+1}}(q;q)_{N_{i+1}-N_{i+2}}\cdots(q;q)_{N_{k-2}-N_{k-1}}(q;q)_{N_{k-1}}}.
\end{split}
\end{equation}
Applying Lemma \ref{bl1} to \eqref{BP6} $i-1$ times, we  get the Bailey pair $(\alpha^{(2k-i)}_n(1,q),\beta^{(2k-i)}_n(1,q))$,
\begin{equation}\label{BP0}
\begin{split}
&\alpha^{(2k-i)}_n(1,q)=\left\{
             \begin{array}{ll}
               1, & \hbox{if $n=0$}, \\[5pt]
               (-1)^nq^{\frac{2k-1}{2}n^2}(q^{-\frac{2k-2i+1}{2}n}
               +q^{\frac{2k-2i-1}{2}n})(1+q^n)/2, & \hbox{if $n\geq1$,}
             \end{array}
               \right.\\[10pt]
&\beta^{(2k-i)}_n(1,q)=\sum_{n\geq N_{2}\geq \cdots\geq N_{k-1}\geq 0}\frac{q^{N^2_{2}+N^2_{3}+\cdots+N^2_{k-1}+N_{i+1}+\cdots+N_{k-1}}(1+q^{N_i})}
{2(q;q)_{n-N_{2}}(q;q)_{N_{2}-N_{3}}\cdots(q;q)_{N_{k-2}-N_{k-1}}(q;q)_{N_{k-1}}}.
\end{split}
\end{equation}
Plugging \eqref{BP0} into Lemma \ref{bailey lemma3}, we get a Bailey pair $(\alpha^{(2k-i+1)}_n(1,q),\beta^{(2k-i+1)}_n(1,q))$,
\begin{equation*}
\begin{split}
&\alpha^{(2k-i+1)}_n(1,q)=\left\{
             \begin{array}{ll}
               1, & \hbox{if $n=0$}, \\[5pt]
               (-1)^nq^{(2k-1)n^2}(q^{-2(k-i)n}+q^{2(k-i)n}), & \hbox{if $n\geq1$,}
             \end{array}
               \right.\\[10pt]
&\beta^{(2k-i+1)}_n(1,q)=\sum_{n\geq N_{1}\geq \cdots\geq N_{k-1}\geq 0}\frac{(-q;q)_{2N_1-1}q^{N_1+2(N^2_{2}+N^2_{3}+\cdots
+N^2_{k-1}+N_{i+1}+\cdots+N_{k-1})}(1+q^{2N_i})}
{(q^2;q^2)_{n-N_1}(q^2;q^2)_{N_1-N_{2}}\cdots
(q^2;q^2)_{N_{k-2}-N_{k-1}}(q^2;q^2)_{N_{k-1}}}.
\end{split}
\end{equation*}
By the definition of a Bailey pair, letting $n\rightarrow \infty$ and multiplying both sides by $(q^2;q^2)_\infty$, we obtain
\begin{equation}\label{the18-temp}
\begin{split}
&\sum_{N_{1}\geq \cdots\geq N_{k-1}\geq 0}\frac{(-q;q)_{2N_1-1}q^{N_1+2(N^2_{2}+N^2_{3}+\cdots+N^2_{k-1}+N_{i+1}+\cdots+N_{k-1})}(1+q^{2N_i})}
{(q^2;q^2)_{N_1-N_{2}}(q^2;q^2)_{N_{2}-N_{3}}
\cdots(q^2;q^2)_{N_{k-2}-N_{k-1}}(q^2;q^2)_{N_{k-1}}}\\
&=\frac{(-q;q)_\infty}{(q;q)_\infty}\left(1+\sum_{n=1}^\infty
(-1)^nq^{(2k-1)n^2}(q^{-2(k-i)n}+q^{2(k-i)n})\right).
\end{split}
\end{equation}
Using Jacobi's triple product identity, we  see that
\begin{align}\label{the18-temp2}
&1+\sum_{n=1}^\infty
(-1)^nq^{(2k-1)n^2}(q^{-2(k-i)n}+q^{2(k-i)n}) =(q^{2i-1},q^{4k-2i-1},q^{4k-2};q^{4k-2})_\infty.
\end{align}
Submitting \eqref{the18-temp2} into \eqref{the18-temp}, and noting that
\[(-q;q)_{2N_1-1}q^{N_1}=(-q^{2-2N_1};q^2)_{N_1-1}
(-q^{1-2N_1};q^2)_{N_1}q^{2N^2_1},\]
we obtain \eqref{over3}. Thus we complete the proof of Theorem \ref{Gordon-Gollnitz-odd-overp-eqn}. \qed

\section{The G\"ollnitz-Gordon marking of an overpartition}

In this section, we  first introduce the G\"ollnitz-Gordon marking of an overpartition and then give an outline of the proof of Theorem \ref{Gollnitz-odd-1}.

Kur\c{s}ung\"{o}z \cite{Kursungoz-2010} introduced the notion of the Gordon marking of an ordinary partition and gave a combinatorial proof of Theorem \ref{R-R-G-left-e}.  The Gordon marking of an ordinary partition $\eta$ is an assignment of positive integers (marks) to the parts of  $\eta$ from  smallest to  largest such that the marks are as small as possible subject to the condition that equal or consecutive parts are assigned different marks. More precisely,    let $\eta=(\eta_1,\eta_2,\ldots,\eta_\ell)$ where $1\leq \eta_1\leq \eta_2 \leq \cdots \leq  \eta_\ell$.   Assign $1$ to $\eta_1$. For $p>1$,  assume that $s$ is the least positive integer that is not used to mark the parts   $\eta_j$ with $\eta_p-\eta_j\leq 1$  for $j<p$. Then, we assign   $s$ to $\eta_p$.  For example, the Gordon marking of $\eta=(1,1,2,2,2,3,4,5,5,6,6,6)$ is
\[G(\eta)=(1_1,1_2,2_3,2_4,2_5,3_1,4_2,5_1,5_3,6_2,6_4,6_5).\]
The Gordon marking of an ordinary partition can also  be represented by an array  where the column indicates the value of parts and the row (counted from bottom to top) indicates the mark
listed outside the brackets, so the Gordon marking of $\eta$ above can be represented as:
\[G(\eta)=
\begin{array}{ccc}\left[
\begin{array}{cccccccc}
&2&&&&6\\
&2&&&&6\\
&2&&&5&\\
1&&&4&&6\\
1&&3&&5&
\end{array}
\right]&\begin{array}{c}5\\4\\3\\2\\1\end{array}.\end{array}
\]
We call this array the Gordon marking representation of an ordinary partition. Let $N_r$ denote the number of parts in the $r$-th row of the Gordon marking representation of an ordinary partition $\eta$, that is, the number of $r$-marked parts in the Gordon marking of  $\eta$. From the definition of the Gordon marking of an ordinary partition, it is not difficult  to see that $N_1\geq N_2 \geq \cdots$. Furthermore, if $\eta$ is counted by $B_{k,i}(n)$  in Theorem \ref{R-R-G}, then there are no $k$ or greater marked parts in the Gordon marking of $\eta$, namely, there are at most $k-1$ rows in  the Gordon marking representation of $\eta$.

 Let $\mathbb{B}_{N_1,\ldots,N_{k-1};i}(n)$ denote the set of partitions $\eta$ counted by $B_{k,i}(n)$   such that there are $N_r$ parts in the $r$-th row of the Gordon marking representation of $\eta$ for $1\leq r\leq k-1$. Define \[\mathbb{B}_{N_1,\ldots,N_{k-1};i}=\bigcup_{n\geq0}
\mathbb{B}_{N_1,\ldots,N_{k-1};i}(n).\]

 Kur\c{s}ung\"{o}z \cite{Kursungoz-2010} established the following identity by introducing backward and forward moves defined on the Gordon marking of an ordinary partition.
\begin{equation}\label{GENREATING-EEK}
\sum_{\eta\in\mathbb{B}_{N_1,\ldots,N_{k-1};i}}
q^{|\eta|}=
\frac{
q^{N^2_1+\cdots+N^2_{k-1}+N_{i}+\cdots+N_{k-1}}}
{(q;q)_{N_1-N_2}\cdots(q;q)_{N_{k-2}-N_{k-1}}
(q;q)_{N_{k-1}}},
\end{equation}
where   $|\eta|$ denotes the sum of parts of $\eta$.

Let $\mathbb{B}_{k,i}(m,n)$ denote the set of partitions counted by ${B}_{k,i}(m,n)$ in Theorem \ref{R-R-G-left-e}. From the definition of $\mathbb{B}_{N_1,\ldots,N_{k-1};i}(n)$, we see that

\[\mathbb{B}_{k,i}(m,n)=\bigcup_{\substack{N_1\geq \cdots \geq N_{k-1}\geq 0\\
N_1+\cdots+N_{k-1}=m}}\mathbb{B}_{N_1,\ldots,N_{k-1};i}(n),\]
so
 \begin{equation}\label{GENREATING-tempt}
 \sum_{m,\,n\geq 0}B_{k,i}(m,n)x^mq^n=\sum_{N_1\geq N_2\geq \ldots \geq N_{k-1}\geq0}x^{N_1+N_2+\ldots+N_{k-1}}
 \sum_{\eta\in\mathbb{B}_{N_1,\ldots,N_{k-1};i}}q^{|\eta|}.
\end{equation}
 Therefore, inserting \eqref{GENREATING-EEK} into \eqref{GENREATING-tempt} to give rise to Theorem \ref{R-R-G-left-e}.

To show Theorem \ref{Gollnitz-odd-1}, we introduce   the G\"ollnitz-Gordon marking of an overpartition which is different from the Gordon marking of an overpartition introduced by Chen, Sang and Shi \cite{Chen-Sang-Shi-2013}. It should be mentioned that an ordinary partition can be marked with   G\"ollnitz-Gordon marking,  but the G\"ollnitz-Gordon marking of an ordinary partition is different from the   Gordon marking of an ordinary partition.

In the remainder of this paper,  we write an overpartition $\lambda$ as the form $(\lambda_1,\,\lambda_2,\ldots, \lambda_\ell)$ where $\lambda_1\leq \lambda_2\leq \cdots \leq \lambda_\ell$ are ranked  in the following order,
\begin{equation}\label{order}
\overline{1}<1<\overline{2}<2<\cdots.
\end{equation}

The $\lambda_j$  is called  the $j$-th part of an overpartition $\lambda$. Denote the size of  $\lambda_j$ by  $|\lambda_j|$. If $|\lambda_j|=a_j$, then we write $\lambda_j=\overline{a_j}$ to indicate that $\lambda_j$ is an overlined part and write $\lambda_j=a_j$ to indicate that $\lambda_j$ is a non-overlined part.

\begin{defi}[G\"ollnitz-Gordon marking]\label{pf-1} The G\"ollnitz-Gordon marking of an overpartition $\lambda$ is an assignment of positive integers (marks) to parts of $\lambda=(\lambda_1,\lambda_2,\ldots,\lambda_\ell)$ from  smallest to  largest. Assign $1$ to $\lambda_1$. For $p>1$, assume that the part $\lambda_j$  has been assigned  a mark for $j<p$.  We consider the following four cases:

$(1)$ If $\lambda_p$ is a non-overlined  odd part, then   assign $1$ to $\lambda_p$.

$(2)$ If $\lambda_p$ is an overlined  even part, then  assign $1$ to $\lambda_p$.

$(3)$ If $\lambda_p$ is an overlined   odd part, and let $s$ be the least positive integer that is not used to mark the parts $\lambda_j$ with $|\lambda_p|-|\lambda_j|=1$ for $j<p$,  then assign $s$ to $\lambda_p$.

$(4)$ If $\lambda_p$ is a non-overlined  even part, say  $\lambda_p=2t+2$, and define

 $\diamondsuit$   $f$ to be the least positive integer that is not used to mark the parts $\lambda_j$  with $|\lambda_p|-|\lambda_j|\leq 2$  for $j<p$,

 $\diamondsuit$   $g$  to be the least positive integer that has been used to mark the parts   $\lambda_j$ with $|\lambda_p|-|\lambda_j|=2$    for $j<p$. If such $\lambda_j$  does not occur in $\lambda$,   then set $g=0$,

then we may assign $f$ or $g$ to $\lambda_p$ by considering the following two subcases:

 $(4.1)$ If $\lambda$ satisfies  four conditions simultaneously: $(i)$ $g\geq 2$;  $(ii)$ the mark of $\lambda_{p-1}$ is $g-1$;  $(iii)$   $2t+1$ or   $\overline{2t+2}$ occurs in  $\lambda$;  $(iv)$   $\overline{2t+1}$ does not occur in $\lambda$, then  assign $g$ to $\lambda_p$;

  $(4.2)$ Otherwise,    assign $f$ to $\lambda_p$.
\end{defi}

For example, we consider the overpartition
  \[ \begin{array}{lllllllllllllllll}
  &\lambda_1,& \lambda_2,& \lambda_3,& \lambda_4,& \lambda_5,& \lambda_6,& \lambda_7,& \lambda_8,& \lambda_9,& \lambda_{10}\\[5pt]
  &\downarrow&\downarrow&\downarrow&\downarrow&\downarrow&\downarrow
  &\downarrow&\downarrow&\downarrow&\downarrow\\[5pt]
 \lambda=(&1,&1,&\overline{2},&2,&\overline{3},
&\overline{4},&6,&7,&8,&8&).
\end{array}\]

By Definition \ref{pf-1}, we see that   $\lambda_1=1$, $\lambda_2=1$ and $\lambda_8=7$ should be marked with 1 since they are non-overlined  odd parts. On the other hand,   $\lambda_3=\overline{2}$ and $\lambda_6=\overline{4}$ are overlined even parts, so they are also marked with 1. Hence, we have   \[ \begin{array}{lllllllllllllllll}
  &\lambda_1,& \lambda_2,& \lambda_3,& \lambda_4,& \lambda_5,& \lambda_6,& \lambda_7,& \lambda_8,& \lambda_9,& \lambda_{10}\\[5pt]
  &\downarrow&\downarrow&\downarrow&\downarrow&\downarrow&\downarrow
  &\downarrow&\downarrow&\downarrow&\downarrow\\[5pt]
 (&1_1,&1_1,&\overline{2}_1,&2,&\overline{3},
&\overline{4}_1,&6,&7_1,&8,&8&).
\end{array}\]
The part $\lambda_4=2$ should be marked with $2$ since it is a non-overlined  even part and it does not satisfy the conditions   (4.1) in Definition \ref{pf-1}. The part $\lambda_5=\overline{3}$ is marked with $3$ since it is an overlined odd part and $\lambda_3=\overline{2}$ and $\lambda_4=2$ are marked with 1 and 2 respectively. The part $\lambda_7=6$ is marked with $2$ since   it is a non-overlined  even part and it does not satisfy the conditions   (4.1) in Definition \ref{pf-1}. The part $\lambda_9=8$ is marked with $2$ since   it is a non-overlined even part and it satisfies the conditions   (4.1) in Definition \ref{pf-1}. The part $\lambda_{10}=8$ is marked with $3$ since   it is a non-overlined even part and it does not  satisfy the conditions   (4.1) in Definition \ref{pf-1}. So the G\"ollnitz-Gordon marking of $\lambda$  is
\[GG(\lambda)=(1_1,1_1,{\overline{2}}_1,2_2,{\overline{3}}_3,
\overline{4}_1,6_2,7_1,8_2,8_3).\]
 It can also be represented   by an array where column indicates the size of parts, and the row (counted from bottom to top) indicates the mark listed outside the brackets, so the G\"ollnitz-Gordon marking of $\lambda$ above would be
\begin{equation}\label{example1}
GG(\lambda)=\begin{array}{ccc}\left[
\begin{array}{cccccccccccccccccccccc}
 &&\overline{3}&&&&&8&\\
 &2&&&&6&&8&\\
 1^2&\overline{2}&&\overline{4}&&&7&
\end{array}
\right]&\begin{array}{c}3\\2\\1\end{array}.\end{array}
\end{equation}

Similarly, we will call this diagram the G\"ollnitz-Gordon marking representation of an overpartition.  Note that   non-overlined odd parts would be repeated in the first row of
  the \GG \  representation of $\lambda$, so for $t\geq 2$,  we  will use ${(2j+1)}^t$  to  denote that there are $t$  $(2j+1)$'s  in the first row of the \GG \ representation of $\lambda$, and ${\overline{2j+1}}^t$ to denote that there are  a $\overline{2j+1}$ and $t-1$    $(2j+1)$'s  in the first row of the \GG \ representation of $\lambda$.

Let $N_r$ denote the number of  parts in the $r$-th row of  the G\"ollnitz-Gordon marking representation of an overpartition. From the definition of   G\"ollnitz-Gordon marking, it is not hard  to show that $N_1\geq N_2\geq\cdots$.  For the example above, we have $N_1=5,\ N_2=3,$ and $N_3=2.$

From the definition of \GG, we see that if $\lambda$ is counted by $O_{k,i}(n)$ in Theorem \ref{Gordon-Gollnitz-odd-overp}, then  $f_{\overline{1}}(\lambda)+f_2(\lambda)\leq i-1$ and  there are no $k$ or greater marked parts in the
the G\"ollnitz-Gordon marking of $\lambda$, that is, there are at most $k-1$ rows in the   G\"ollnitz-Gordon marking representation of $\lambda$, and vice versa.
More precisely, we have following proposition.

\begin{pro}For $k\geq i\geq1$, an overpartition $\lambda$ is counted by $ {{O}}_{k,i}(n)$  if and only if  the number of occurrences  of $\bar{1}$ and $2$ in $\lambda$ is not exceed $i-1$ and there are at most $k-1$ rows in the  G\"ollnitz-Gordon marking representation of $\lambda$.
\end{pro}

It should be noted that for the parts  $2t+2$ in $\lambda$, if  $f_{\overline{2t+1}}(\lambda)=0$, $f_{\overline{2t}}(\lambda)=0$ and
\[ f_{2t}(\lambda)+f_{2t+2}(\lambda)= k-1,\]
and the least positive  integer $q$  that has been used to mark the parts  $2t$ in $\lambda$ is greater than $1$, and there is at least one $2t+1$ or $\overline{2t+2}$   in $\lambda$ which will be marked with 1,  then  the marks of $2t+2$ will be less than $k$ since there is a   $2t+2$ in $\lambda$ marked with $g$ which satisfies  the conditions  in (4.1) of Definition \ref{pf-1}. This is the reason that the marking of non-overlined even parts in the  definition of \GG \   is more complicated.

For $k\geq i\geq 1$, let  $\mathbb{O}_{k,i}(m,n)$ denote the set of  overpartitions counted by $O_{k,i}(m,n)$,   we will classify $\mathbb{O}_{k,i}(m,n)$  by considering whether the smallest part of an overpartition in $\mathbb{O}_{k,i}(m,n)$ is a non-overlined odd part or an overlined even part. Note that the parts of an overpartition are ordered by \eqref{order}. Let $\mathbb{F}_{k,i}(m,n)$ denote the set of overpartitions in $\mathbb{O}_{k,i}(m,n)$ for which the smallest part is an overlined odd part or a non-overlined even part, and let $\mathbb{H}_{k,i}(m,n)$ denote the set of overpartitions in $\mathbb{O}_{k,i}(m,n)$ with the smallest part being a non-overlined odd part or an overlined even part. Obviously,
\begin{equation}\label{2N-1}
\mathbb{O}_{k,i}(m,n)=\mathbb{F}_{k,i}(m,n)\cup\mathbb{H}_{k,i}(m,n).
\end{equation}

Let ${F}_{k,i}(m,n)=|\mathbb{F}_{k,i}(m,n)|$ and ${H}_{k,i}(m,n)=|\mathbb{H}_{k,i}(m,n)|$. Then
\begin{equation}\label{2N-2}
{O}_{k,i}(m,n)={F}_{k,i}(m,n)+{H}_{k,i}(m,n).
\end{equation}

There is a relation between ${F}_{k,i}(m,n)$ and  ${H}_{k,i}(m,n)$.

\begin{lem}\label{2N1}
For $k\geq i\geq2$,
\begin{equation}\label{2N-3}
{F}_{k,i}(m,n)={H}_{k,i-1}(m,n).
\end{equation}
For $k\geq 1$,
\begin{equation}\label{2N-4}
{F}_{k,1}(m,n)={H}_{k,k}(m,n-2m).
\end{equation}
\end{lem}

\pf For $k\geq i\geq2$, there is a simple bijection between $\mathbb{F}_{k,i}(m,n)$ and $\mathbb{H}_{k,i-1}(m,n)$. Let $\sigma$ be an overpartition in $\mathbb{F}_{k,i}(m,n)$, we consider two cases: If the smallest part of $\sigma$ is  an overlined odd part, say $\overline{2t+1}$, then  change it to a non-overlined odd part ${2t+1}$.  If the smallest part of $\sigma$ is a non-overlined even part, say ${2t}$, then change the first ${2t}$ of $\sigma$   to an overlined even part $\overline{2t}$. In either case, we   obtain an overpartition $\pi$ in $\mathbb{H}_{k,i-1}(m,n)$. Furthermore, it is easy to see that this process is reversible, and so this map is a bijection. Hence \eqref{2N-3} holds for $k\geq i\geq2$.

For $k\geq  1$, we will give a bijection between $\mathbb{F}_{k,1}(m,n)$ and $\mathbb{H}_{k,k}(m,n-2m)$. For an overpartition $\sigma\in\mathbb{F}_{k,1}(m,n)$, by the definition of $\mathbb{F}_{k,1}(m,n)$, we see that $\sigma$ has $m$ parts and the size  of each part of $\sigma$ is greater than $2$. There are two cases: If the smallest part of $\sigma$ is  an overlined odd part, say $\overline{2t+1}$, where $t\geq 1$, then change it to a non-overlined odd part ${2t+1}$. If the smallest part of $\sigma$ is   a non-overlined even part,  then change one of the smallest parts, say  ${2t}$, where $t\geq 2$ to an overlined even part $\overline{2t}$. In either case, we obtained a new overpartition $\rho$ for which the size of each part  is greater than $2$. We then  subtract $2$ from each part of $\rho$ to obtain the resulting overpartition $\pi$ in $\mathbb{H}_{k,k}(m,n-2m)$. It is evident to see that this process is reversible, and so this map is a bijection. Hence   we arrive at \eqref{2N-4}. This completes the proof. \qed

Using the relation \eqref{2N-2},  it is easy to find that the generating function of $O_{k,i}(m,n)$ can be deduced from the generating functions of $F_{k,i}(m,n)$ and $H_{k,i}(m,n)$. In light of Lemma \ref{2N1}, we see that the generating function of $H_{k,i}(m,n)$ can be obtained from the generating function of $F_{k,i}(m,n)$. Hence,  it suffices to derive the following generating  function of $F_{k,i}(m,n)$ in order to prove Theorem \ref{Gollnitz-odd-1}.

\begin{thm}\label{thm-2N-1}
For $k\geq i\geq1,$
\begin{equation}\label{thm-2N1}
\begin{split}
& \displaystyle\sum_{m,n\geq0}F_{k,i}(m,n)x^mq^n\\
&=\sum_{N_{1}\geq \cdots\geq N_{k-1}\geq 0}\frac{(-q^{2-2N_1};q^2)_{N_1-1}(-q^{1-2N_1};q^2)_{N_1}
q^{2(N^2_1+\cdots+N^2_{k-1}+N_{i}+\cdots+N_{k-1})}
x^{N_1+\cdots+N_{k-1}}}
{(q^2;q^2)_{N_1-N_2}\cdots(q^2;q^2)_{N_{k-2}-N_{k-1}}
(q^{2};q^{2})_{N_{k-1}}}.
\end{split}
\end{equation}
\end{thm}

In this paper, we will give a combinatorial proof of Theorem \ref{thm-2N-1} based on the \GG \  of an overpartition. Let   $\mathbb{F}_{N_1,\ldots,N_{k-1};i}(n)$ denote the set of overpartitions  $\lambda$ in $\mathbb{F}_{k,i}(m,n)$ such that there are $N_r$ $r$-marked parts in the G\"ollnitz-Gordon marking of $\lambda$ for $1\leq r\leq k-1$.

Set
\[\mathbb{F}_{N_1,\ldots,N_{k-1};i}
=\bigcup_{n\geq0}\mathbb{F}_{N_1,\ldots,N_{k-1};i}(n).\]

From the definition of $\mathbb{F}_{N_1,\ldots,N_{k-1};i}$, it is evident to see that

\[\mathbb{F}_{k,i}(m,n)=\bigcup_{\substack{N_1\geq \cdots \geq N_{k-1}\geq 0\\
N_1+\cdots+N_{k-1}=m}}\mathbb{F}_{N_1,\ldots,N_{k-1};i}(n).\]

This leads to
 \begin{equation}\label{GENREATING-tempt-f}
 \sum_{m,\,n\geq 0}F_{k,i}(m,n)x^mq^n=\sum_{N_1\geq N_2\geq \ldots \geq N_{k-1}\geq0}x^{N_1+\cdots+N_{k-1}}\sum_{\lambda\in\mathbb{F}_{N_1,\ldots,N_{k-1};i}}q^{|\lambda|}.
\end{equation}

Hence Theorem \ref{thm-2N-1} immediately follows when we show that for $k\geq i\geq 1$ and  $N_1\geq N_2\geq \cdots \geq N_{k-1}\geq 0$,
\begin{align}\label{N1e}
&\sum_{\lambda\in\mathbb{F}_{N_1,\ldots,N_{k-1};i}}
q^{|\lambda|}\nonumber\\[3pt]
&\quad =\frac{(-q^{2-2N_1};q^2)_{N_1-1}(-q^{1-2N_1};q^2)_{N_1}
q^{2(N^2_1+\cdots+N^2_{k-1}+N_{i}+\cdots+N_{k-1})}
}
{(q^2;q^2)_{N_1-N_2}\cdots(q^2;q^2)_{N_{k-2}-N_{k-1}}
(q^{2};q^{2})_{N_{k-1}}}.
\end{align}

 It turns out that the proof of \eqref{N1e} is more complicated than the proof of \eqref{GENREATING-EEK} due to Kur\c{s}ung\"{o}z. To prove \eqref{N1e}, we require to build  two bijections. More precisely,
let $\mathbb{G}_{N_1,\ldots,N_{k-1};i}(n)$ denote the set of overpartitions in $\mathbb{F}_{N_1,\ldots,N_{k-1};i}(n)$  for which there are no   overlined even parts and non-overlined odd parts. Let $\mathbb{E}_{N_1,\ldots,N_{k-1};i}(n)$ denote the set of
overpartitions  in $\mathbb{G}_{N_1,\ldots,N_{k-1};i}(n)$ for which there are no  overlined odd parts.

Set
\[\mathbb{G}_{N_1,\ldots,N_{k-1};i}
=\bigcup_{n\geq0}\mathbb{G}_{N_1,\ldots,N_{k-1};i}(n),\]
and
\[\mathbb{E}_{N_1,\ldots,N_{k-1};i}=\bigcup_{n\geq 0}\mathbb{E}_{N_1,\ldots,N_{k-1};i}(n).\]

We will give bijective proofs of the following two relations in  Section 4 and Section 5 respectively.

\begin{lem}\label{N-1}
For $k\geq i\geq1$ and $N_1\geq N_2\geq \cdots \geq N_{k-1}\geq 0$,
\begin{equation}\label{N1}
\sum_{\lambda\in\mathbb{F}_{N_1,\ldots,N_{k-1};i}}
q^{|\lambda|}=(-q^{2-2N_1};q^2)_{N_1-1}\sum_{\mu\in
\mathbb{G}_{N_1,\ldots,N_{k-1};i}}
q^{|\mu|}.
\end{equation}

\end{lem}

 \begin{lem}\label{N-2}
For $k\geq i\geq1$ and $N_1\geq N_2\geq \cdots \geq N_{k-1}\geq 0$,
\begin{equation}\label{N2e}
\sum_{\mu\in\mathbb{G}_{N_1,\ldots,N_{k-1};i}}
q^{|\mu|}=(-q^{1-2N_1};q^2)_{N_1}\sum_{\nu\in\mathbb{E}_{N_1,\ldots,N_{k-1};i}}
q^{|\nu|}.
\end{equation}
\end{lem}

In Section 6, we first give a proof of Theorem \ref{thm-2N-1}  by using Lemma \ref{N-1} and Lemma \ref{N-2}, as well as  \eqref{GENREATING-EEK} due to Kur\c{s}ung\"{o}z. Then we complete the proof of  Theorem \ref{Gollnitz-odd-1} by using Theorem \ref{thm-2N-1}, together with  Lemma \ref{2N1} and the relation \eqref{2N-2}. In the remaining part of this paper, we mark  parts of an  overpartition in \GG.

\section{Proof of Lemma \ref{N-1}}

Let $\mathbb{R}_N$ denote the set of partitions $\tau=(\tau_1,\tau_2,\ldots, \tau_\ell)$ with distinct negative even parts which lay in $[-2N,-2]$, that is, $\tau_j$ is negative and even for $1\leq  j\leq \ell$ and $-2N\leq \tau_1<\tau_2<\cdots<\tau_\ell\leq -2$. It is easy to see that the generating function of partitions in  $ \mathbb{R}_N$:
\[\sum_{\tau \in \mathbb{R}_{N}}q^{|\tau|}=(1+q^{-2N})(1+q^{2-2N})\cdots (1+q^{-2})=(-q^{-2N};q^2)_{N}.
\]
Thus, Lemma \ref{N-1} is equivalent to the following combinatorial statement.

\begin{thm}\label{N-1-COM} For $k\geq i\geq 1$ and $N_1\geq N_2\geq\cdots \geq N_{k-1}\geq 0$, there is  a bijection $\Phi$ between $\mathbb{{F}}_{N_1,\ldots,N_{k-1};i}$ and $\mathbb{R}_{N_1-1}\times\mathbb{{G}}_{N_1,\ldots,N_{k-1};i}$ such that   for $\lambda \in \mathbb{F}_{N_1,\ldots,N_{k-1};i}$ and $\Phi(\lambda)=(\tau,\mu)\in \mathbb{R}_{N_1-1}\times\mathbb{{G}}_{N_1,\ldots,N_{k-1};i}$, we have $|\lambda|=|\tau|+|\mu|.$
\end{thm}

   Observe that $\mathbb{G}_{N_1,\ldots,N_{k-1};i}$ is the set of overpartitions in $\mathbb{F}_{N_1,\ldots,N_{k-1};i}$  for which there are no   overlined even parts and non-overlined odd parts, so the key point  in the construction of the bijection $\Phi$  is   to  remove  overlined even parts and non-overlined odd parts from an overpartition in $\mathbb{F}_{N_1,\ldots,N_{k-1};i}$ to obtain a new overpartition in $\mathbb{G}_{N_1,\ldots,N_{k-1};i}$. To this end, we will first define three subsets $\mathbb{F}_{N_1,\ldots,N_{k-1};i,p}$, $\overline{\mathbb{F}}_{N_1,\ldots,N_{k-1};i,p}$ and $\overrightarrow{\mathbb{F}}_{N_1,\ldots,N_{k-1};i,p}$ of $\mathbb{F}_{N_1,\ldots,N_{k-1};i}$. Then we build a bijection $\Phi_p$ between $\mathbb{F}_{N_1,\ldots,N_{k-1};i,p}$ and  $\overline{\mathbb{F}}_{N_1,\ldots,N_{k-1};i,p}$ and a bijection $\Phi_{(p)}$ between $\mathbb{F}_{N_1,\ldots,N_{k-1};i,p}$ and  $\overrightarrow{\mathbb{F}}_{N_1,\ldots,N_{k-1};i,p}$. It turns out that    $\Phi_{(p)}$ can be obtained by iteratively using the bijection $\Phi_p$ and plays a crucial role in the construction of the bijection $\Phi$ in Theorem \ref{N-1-COM}.

Let $\lambda$ be an overpartition in $\mathbb{F}_{N_1,\ldots,N_{k-1};i}$. For $1\leq r\leq k-1$, define $\lambda^{(r)}=(\lambda^{(r)}_1,\,\lambda^{(r)}_2,\ldots,\break \lambda^{(r)}_{N_r})$ to be the $r$-th sub-overpartition of $\lambda$ whose parts are $r$-marked parts in the G\"ollnitz-Gordon marking of $\lambda$, where $\lambda^{(r)}_1\leq \lambda^{(r)}_2\leq \cdots\leq \lambda^{(r)}_{N_r}$. The $\lambda^{(r)}_j$ are called the $j$-th part of the   $r$-th sub-overpartition $\lambda^{(r)}$ of $\lambda$.

Let $k\geq i\geq 1$ and $N_1\geq N_2\geq\cdots \geq N_{k-1}\geq 0$ be given. For $1<p \leq N_1$, the subsets $\mathbb{F}_{N_1,\ldots,N_{k-1};i,p}$, $\overline{\mathbb{F}}_{N_1,\ldots,N_{k-1};i,p}$ and $\overrightarrow{\mathbb{F}}_{N_1,\ldots,N_{k-1};i,p}$ are described by using the first sub-overpartition $\lambda^{(1)}=(\lambda^{(1)}_1,\,\lambda^{(1)}_2,\ldots, \lambda^{(1)}_{N_1})$ of an overpartition $\lambda$ in $\mathbb{F}_{N_1,\ldots,N_{k-1};i}$, where $\lambda^{(1)}_1\leq \lambda^{(1)}_2\leq \cdots\leq \lambda^{(1)}_{N_1}$.

 $\lozenge$ $\mathbb{F}_{N_1,\ldots,N_{k-1};i,p}$ is the set of   overpartitions $\lambda$ in $\mathbb{F}_{N_1,\ldots,N_{k-1};i}$ such  that (1) $\lambda_p^{(1)}$  is a non-overlined odd part or an overlined even part; (2) $\lambda_j^{(1)}$  is an overlined odd part or a non-overlined even part, where  $ p+1\leq j\leq N_1$.

$\lozenge$   $\mathbb{\overline{F}}_{N_1,\ldots,N_{k-1};i,p}$ is the set of   overpartitions $\lambda$ in $\mathbb{F}_{N_1,\ldots,N_{k-1};i}$ such that
    (1) $\lambda_p^{(1)}$  is an overlined odd part or a non-overlined even part; (2) $\lambda_{p+1}^{(1)}$  is a non-overlined odd part or an overlined even part; (3)  $\lambda_j^{(1)}$ is an overlined odd part or a non-overlined even part, where $p+2\leq j\leq N_1$.

$\lozenge$   $\overrightarrow{\mathbb{F}}_{N_1,\ldots,N_{k-1};i,p}$ is the set of   overpartitions $\lambda$ in $\mathbb{F}_{N_1,\ldots,N_{k-1};i}$ such that  $\lambda_j^{(1)}$ is an overlined odd part or a non-overlined even part, where  $p\leq j\leq N_1$.

By definition, it is easy to see that for $1\leq p\leq N_1-2,$
 \[\mathbb{\overline{F}}_{N_1,\ldots,N_{k-1};i,p}  \subseteq \mathbb{F}_{N_1,\ldots,N_{k-1};i,{p+1}} \subseteq  \overrightarrow{\mathbb{F}}_{N_1,\ldots,N_{k-1};i,p+2}.\]

We are ready to present the bijection $\Phi_p$ between   $\mathbb{F}_{N_1,\ldots,N_{k-1};i,p}$ and  $\overline{\mathbb{F}}_{N_1,\ldots,N_{k-1};i,p}$ and the bijection $\Phi_{(p)}$ between $\mathbb{F}_{N_1,\ldots,N_{k-1};i,p}$ and  $\overrightarrow{\mathbb{F}}_{N_1,\ldots,N_{k-1};i,p}$.  The following lemma  gives the bijection $\Phi_p$, which will be proved at the end of this section.

 \begin{lem}\label{N-1-lem}
For $1< p\leq N_1$, there is a bijection $\Phi_p$ between   $\mathbb{F}_{N_1,\ldots,N_{k-1};i,p}$ and $\mathbb{\overline{F}}_{N_1,\ldots,N_{k-1};i,p}$. Furthermore, for $\lambda \in \mathbb{{F}}_{N_1,\ldots,N_{k-1};i,p}$ and $\mu=\Phi_p(\lambda)\in \mathbb{\overline{F}}_{N_1,\ldots,N_{k-1};i,p}$, we have
\begin{equation}
|\mu|=|\lambda|+2,\quad \text{and} \quad \lambda_j^{(1)}=\mu_j^{(1)} \quad \text{for } j\neq p,p+1,
\end{equation}
where $\lambda_j^{(1)}$ {\rm(}resp. $\mu_j^{(1)}${\rm )} is the $j$-th part of the first sub-overpartition of $\lambda$ {\rm(}resp. $\mu${\rm )}.

 \end{lem}

Taking the above bijection $\Phi_p$  in successive to give the following lemma.

 \begin{lem}\label{N-2-lem}
For $1< p\leq N_1$, there is a bijection $\Phi_{(p)}$ between  $\mathbb{F}_{N_1,\ldots,N_{k-1};i,p}$ and  $\overrightarrow{\mathbb{F}}_{N_1,\ldots,N_{k-1};i,p}$. Furthermore, for $\lambda \in \mathbb{{F}}_{N_1,\ldots,N_{k-1};i,p}$ and $\mu=\Phi_p(\lambda)\in \overrightarrow{\mathbb{F}}_{N_1,\ldots,N_{k-1};i,p}$, we have
\begin{equation}\label{N-2-lem-eq}
|\mu|=|\lambda|+2N_1-2p+2,\quad \text{and} \quad \lambda_j^{(1)}=\mu_j^{(1)} \quad \text{for } j<p.
\end{equation}
 \end{lem}

  \pf Define $\Phi_{(p)}=\Phi_{N_1}\Phi_{N_1-1}\cdots \Phi_p$, by Lemma \ref{N-1-lem}, it is easy to verify that $\Phi_{(p)}$ is a bijection between   $\mathbb{F}_{N_1,\ldots,N_{k-1};i,p} $ and  $\overrightarrow{\mathbb{F}}_{N_1,\ldots,N_{k-1};i,p} $ as desired. \qed

Before giving a proof of Lemma \ref{N-1-lem}, we   give a proof of Theorem \ref{N-1-COM} by successively using the bijection $\Phi_{(p)}$ in Lemma \ref{N-2-lem}. \vskip 0.2cm

\noindent{\bf Proof of Theorem \ref{N-1-COM}.}   Let
$\lambda$ be  an overpartition  in $\mathbb{{F}}_{N_1,\ldots,N_{k-1};i}$.  We  aim to define $\Phi(\lambda)=(\tau,\mu)$ such that  $\tau$ is a partition in $\mathbb{R}_{N_1-1}$ and $\mu$ is an overpartition in $\mathbb{{G}}_{N_1,\ldots,N_{k-1};i}$ satisfying $|\lambda|=|\tau|+|\mu|$. We consider  two cases.

Case 1. If there are no   overlined even parts and  non-overlined odd parts in $\lambda$, then set $\mu=\lambda$ and $\tau=\emptyset$.  It is easy to see that $\mu \in \mathbb{{G}}_{N_1,\ldots,N_{k-1};i}$ and $|\lambda|=|\mu|$.

Case 2. If there are $s\geq 1$ overlined even parts or non-overlined odd parts in $\lambda$, then by the definition of  G\"ollnitz-Gordon marking, these parts are marked with 1. If we assume that $\lambda^{(1)}=(\lambda_1^{(1)},\lambda_2^{(1)},\ldots,\lambda_{N_1}^{(1)})$ is the first sub-overpartition of $\lambda$, then  there are $s$ overlined even parts or non-overlined odd parts in $\lambda^{(1)}$, which are  $\lambda_{j_1}^{(1)}$, $\lambda_{j_2}^{(1)}$, \ldots, $\lambda_{j_{s-1}}^{(1)}$ and $\lambda_{j_s}^{(1)}$, where $1\leq j_1<j_2< \cdots < j_s \leq N_1$. Under this assumption, we see that $\lambda \in \mathbb{{F}}_{N_1,\ldots,N_{k-1};i, j_{s}}$. Note that the smallest part of $\lambda$ is an overlined odd part or a non-overlined even part, so $j_1>1$. Set
  \[\tau=(-2(N_1-j_1+1), -2(N_1-j_2+1),\ldots,\ -2(N_1-j_s+1)),\]
   obviously, $\tau$ is a partition in $\mathbb{R}_{N_1-1}$. The partition $\mu$ can be obtained from $\lambda$ by employing the bijection in Lemma \ref{N-2-lem} $s$ times. We denote the intermediate partitions by $\gamma^{0},\ \gamma^{1},\ldots,\gamma^{s}$ with $\gamma^{0}=\lambda$ and $\gamma^{s}=\mu$.
 For $1\leq b\leq s$,  the intermediate partition $\gamma^{b}$ can be obtained from  $\gamma^{b-1}$ by using $\Phi_{(j_{s-b+1})}$ in Lemma  \ref{N-2-lem}, that is,  for  $1\leq b\leq s$,
 \[\gamma^{b}=\Phi_{(j_{s-b+1})}(\gamma^{b-1}).\]
   Note that $\gamma^{0} \in \mathbb{{F}}_{N_1,\ldots,N_{k-1};i, j_s}$, so  by Lemma   \ref{N-2-lem}, we see that
   \[\gamma^{1} \in \mathbb{{F}}_{N_1,\ldots,N_{k-1};i, j_{s-1}}\text{ and }|\gamma^{1}|=|\lambda|+2(N_1-j_s+1),\]
   and  the first $(j_s-1)$ 1-marked parts in the \GG \  of $\gamma^{1}$ and $\gamma^{0}=\lambda$ are the same.

    Employing Lemma \ref{N-2-lem} repeatedly, we derive that for $1\leq b\leq s-1$,
   \[\gamma^{b} \in \mathbb{{F}}_{N_1,\ldots,N_{k-1};i, j_{s-b}}\text{ and }|\gamma^{b}|=|\lambda|+2\sum_{r=1}^b(N_1-j_{s-r+1}+1),\] and
   \[\gamma^{s} \in \overrightarrow{\mathbb{F}}_{N_1,\ldots,N_{k-1};i, j_{1}}\text{ and }|\gamma^{s}|=|\lambda|+2\sum_{r=1}^s(N_1-j_{s-r+1}+1).\]
   Furthermore, for $1\leq b\leq s$, the first $(j_{s-b+1}-1)$ 1-marked parts in the \GG \  of $\gamma^{b}$ and $\gamma^{0}=\lambda$ are the same. From the preceding assumption, we see that the first $(j_1-1)$ parts in the first sub-overpartition of $\lambda$ are non-overlined even parts or overlined odd parts, and so we derive that   there are no  overlined even parts or non-overlined odd parts in  $\gamma^{s}$. Hence
   \[\mu=\gamma^{s}\in \mathbb{G}_{N_1,\ldots,N_{k-1};i}\text{ and }|\mu|=|\lambda|+2\sum_{r=1}^s(N_1-j_{s-r+1}+1),\] and it is easy to check that $|\tau|+|\mu|=|\lambda|$.  Therefore $\Phi$ is well-defined.

 To prove that  $\Phi$ is a bijection, we shall give a brief description of the inverse map $\Psi$ of $\Phi$.  Let $\mu$ be an overpartition in $\mathbb{G}_{N_1,\ldots,N_{k-1};i}$ and $\tau$ be a partition with distinct negative even parts lying in $[2-2N_1,-2]$. We shall define  $\Psi(\tau,\mu)=\lambda$ such that $\lambda$ is an overpartition in $\mathbb{F}_{N_1,\ldots,N_{k-1};i}$ and $|\lambda|=|\tau|+|\mu|$.  There are two cases.

Case 1. $\tau=\emptyset$. In this event,  set $\lambda=\mu$. Note that $\mathbb{G}_{N_1,\ldots,N_{k-1};i}\subseteq \mathbb{F}_{N_1,\ldots,N_{k-1};i}$, so $\lambda \in \mathbb{F}_{N_1,\ldots,N_{k-1};i}$ and there are no   overlined even parts and non-overlined odd parts in $\lambda$.

Case 2. $\tau\neq \emptyset$. In this event, assume that
\[\tau=(-2(N_1-j_1+1), -2(N_1-j_2+1),\ldots,\ -2(N_1-j_s+1)),\]
 where $1< j_1<j_2<\cdots<j_s\leq N_1.$
  The partition $\lambda$ can be recovered from $\mu$ by  using the bijection  in Lemma  \ref{N-2-lem} $s$ times. We denote the intermediate partitions by $\delta^{0},\,\delta^{1},\ldots,\delta^{s}$ with $\delta^{s}=\mu$ and $\delta^{0}=\lambda$. For $1\leq b\leq s$, the intermediate partition $\delta^{b-1}$ can be obtained from  $\delta^{b}$ by using the   bijection $\Phi^{-1}_{(j_{s-b+1})}$ , that is $\delta^{b-1}=\Phi^{-1}_{(j_{s-b+1})}(\delta^{b})$. Set $\lambda=\delta^0$, it follows from Lemma  \ref{N-2-lem}  that $\lambda$ is  an overpartition in $\mathbb{F}_{N_1,\ldots,N_{k-1};i}$ and $|\lambda|=|\tau|+|\mu|$, and $\Psi(\Phi(\lambda))=\lambda$ for any $\lambda$ in $\mathbb{F}_{N_1,\ldots,N_{k-1};i}$. Hence  $\Phi$ is a bijection between $\mathbb{{F}}_{N_1,\ldots,N_{k-1};i}$ and $\mathbb{R}_{N_1-1}\times\mathbb{{G}}_{N_1,\ldots,N_{k-1};i}$.
 This completes the proof of Theorem \ref{N-1-COM}, and hence Lemma \ref{N-1} is verified.   \qed

It remains to show Lemma \ref{N-1-lem}.  To this end, we shall divide    $\mathbb{F}_{N_1,\ldots,N_{k-1};i,p}$ into four disjoint subsets  $\mathbb{F}^{(l)}_{N_1,\ldots,N_{k-1};i,p}$ ($1\leq l\leq 4$) and divide
  $\mathbb{\overline{F}}_{N_1,\ldots,N_{k-1};i,p}$ into four disjoint subsets  $\mathbb{\overline{F}}^{(l)}_{N_1,\ldots,N_{k-1};i,p}$ ($1\leq l\leq 4$). We then construct the bijection $\Phi_p$ consisting of four bijections $\Phi_{l,p}$ between   $\mathbb{F}^{(l)}_{N_1,\ldots,N_{k-1};i,p}$  and   $\mathbb{\overline{F}}^{(l)}_{N_1,\ldots,N_{k-1};i,p}$, where $1\leq l\leq 4$.

For $1<p\leq  N_1$, let $\lambda^{(1)}=(\lambda^{(1)}_1,\,\lambda^{(1)}_2,\ldots, \lambda^{(1)}_{N_1})$ be the first sub-overpartition of $\lambda$ in $\mathbb{F}_{N_1,\ldots,N_{k-1};i,p}$, by definition,  we see that  $\lambda_p^{(1)}$  is a non-overlined odd part or an overlined even part and $\lambda_j^{(1)}$  is an overlined odd part or a non-overlined even part, where $ p+1\leq j\leq N_1$. The subsets $\mathbb{F}^{(l)}_{N_1,\ldots,N_{k-1};i,p}$ can be described in terms of the first sub-overpartition of $\lambda$.

 \begin{itemize}
  \item[(1)] $\mathbb{F}^{(1)}_{N_1,\ldots,N_{k-1};i,p}$ is the set of overpartitions $\lambda$ in $\mathbb{F}_{N_1,\ldots,N_{k-1};i,p}$ such that (i) $\lambda_p^{(1)}$ is a non-overlined odd part; (ii) there are no non-overlined even parts of size $|\lambda_{p}^{(1)}|+1$ in $\lambda$ when $|\lambda_{p-1}^{(1)}|\leq |\lambda_p^{(1)}|-2$.

 \item[(2)]$\mathbb{F}^{(2)}_{N_1,\ldots,N_{k-1};i,p}$ is the set of overpartitions $\lambda$ in $\mathbb{F}_{N_1,\ldots,N_{k-1};i,p}$ such that (i) $\lambda_p^{(1)}$ is a non-overlined odd part; (ii) there  is at least one non-overlined even part of size $|\lambda_{p}^{(1)}|+1$ in $\lambda$;  (iii) $|\lambda_{p-1}^{(1)}|\leq |\lambda_p^{(1)}|-2$.

 \item[(3)] $\mathbb{F}^{(3)}_{N_1,\ldots,N_{k-1};i,p}$ is the set of overpartitions $\lambda$ in $\mathbb{F}_{N_1,\ldots,N_{k-1};i,p}$ such that (i) $\lambda_p^{(1)}$ is an overlined even part; (ii) there is no overlined odd part of size $|\lambda_p^{(1)}|+1$ in $\lambda$.

 \item[(4)] $\mathbb{F}^{(4)}_{N_1,\ldots,N_{k-1};i,p}$ is the set of overpartitions $\lambda$ in $\mathbb{F}_{N_1,\ldots,N_{k-1};i,p}$ such that (i) $\lambda_p^{(1)}$ is an overlined even part; (ii) there is an overlined odd part of size $|\lambda_p^{(1)}|+1$ in $\lambda$.
 \end{itemize}

For $1<p\leq  N_1$, let $\mu^{(1)}=(\mu^{(1)}_1,\,\mu^{(1)}_2,\ldots, \mu^{(1)}_{N_1})$ be the first sub-overpartition of $\mu$ in $\mathbb{\overline{F}}_{N_1,\ldots,N_{k-1};i,p}$, by definition,  we see that  $\mu_p^{(1)}$  is an overlined odd part or a non-overlined even part, $\mu_{p+1}^{(1)}$  is a non-overlined odd part or an overlined even part, and $\mu_j^{(1)}$ is an overlined odd part or a non-overlined even part, where $p+2\leq j\leq N_1$.  We shall divide  the overpartitions $\mu$ in $\mathbb{\overline{F}}_{N_1,\ldots,N_{k-1};i,p}$ into four disjoint subsets  $\mathbb{\overline{F}}^{(l)}_{N_1,\ldots,N_{k-1};i,p}$, where $1\leq l\leq 4$  in terms of the first sub-overpartition of $\mu$. When $p\geq N_1$, we see that $\mu_{p+1}^{(1)}$ does not occur in $\mu$. For convenience,  set $|\mu_{N_1+1}^{(1)}|=\infty$.

 \begin{itemize}
 \item[(1)] $\mathbb{\overline{F}}^{(1)}_{N_1,\ldots,N_{k-1};i,p}$ is the set of overpartitions $\mu$ in $\mathbb{\overline{F}}_{N_1,\ldots,N_{k-1};i,p}$ such that (i) $\mu_p^{(1)}$ is an overlined odd part;  (ii) there are no non-overlined even parts of size  $|\mu_p^{(1)}|+1$ in $\mu$ when $|\mu_{p+1}^{(1)}|\geq |\mu_p^{(1)}|+2$.

 \item[(2)]$\mathbb{\overline{F}}^{(2)}_{N_1,\ldots,N_{k-1};i,p}$ is the set of overpartitions $\mu$ in $\mathbb{\overline{F}}_{N_1,\ldots,N_{k-1};i,p}$ such that (i) $\mu_p^{(1)}$ is a non-overlined even part; (ii) there is an overlined odd part of size $|\mu_p^{(1)}|+1$ in $\mu$; (iii) there are no non-overlined even parts of size  $|\mu_p^{(1)}|+2$ in $\mu$ when $|\mu_{p+1}^{(1)}|> |\mu_p^{(1)}|+2$.

 \item[(3)] $\mathbb{\overline{F}}^{(3)}_{N_1,\ldots,N_{k-1};i,p}$ is the set of overpartitions $\mu$ in $\mathbb{\overline{F}}_{N_1,\ldots,N_{k-1};i,p}$ such that (i) $\mu_p^{(1)}$ is a non-overlined even part; (ii) there is no overlined odd part of size $|\mu_p^{(1)}|+1$ in $\mu$.

 \item[(4)] $\mathbb{\overline{F}}^{(4)}_{N_1,\ldots,N_{k-1};i,p}$ is the set of overpartitions $\mu$ in $\mathbb{\overline{F}}_{N_1,\ldots,N_{k-1};i,p}$ such that (i) if $\mu_p^{(1)}$ is an overlined odd part, then there is at least one non-overlined even part of size  $|\mu_p^{(1)}|+1$ in $\mu$, and $|\mu_{p+1}^{(1)}|\geq |\mu_p^{(1)}|+2$; (ii) if $\mu_p^{(1)}$ is a non-overlined even part, then  there are an overlined odd part of size $|\mu_p^{(1)}|+1$ and at least one non-overlined even part of size  $|\mu_p^{(1)}|+2$ in $\mu$, and   $|\mu_{p+1}^{(1)}|> |\mu_p^{(1)}|+2$.
 \end{itemize}

 We are now ready to define the bijections $\Phi_{l,p}$ between   $\mathbb{F}^{(l)}_{N_1,\ldots,N_{k-1};i,p}$ and $\mathbb{\overline{F}}^{(l)}_{N_1,\ldots,N_{k-1};i,p}$, where $1\leq l\leq4$.

\begin{lem}\label{theta1}For $1<p\leq N_1$,
there is a bijection $\Phi_{1,p}$ between $\mathbb{F}^{(1)}_{N_1,\ldots,N_{k-1};i,p}$ and  $\mathbb{\overline{F}}^{(1)}_{N_1,\ldots,N_{k-1};i,p}$. Furthermore, for
$\lambda \in \mathbb{{F}}^{(1)}_{N_1,\ldots,N_{k-1};i,p}$  and   $\mu=\Phi_{1,p}(\lambda)\in \overline{\mathbb{F}}^{(1)}_{N_1,\ldots,N_{k-1};i,p},$
 we have
\[|\mu|=|\lambda|+2,\quad \text{and} \quad \lambda_j^{(1)}=\mu_j^{(1)} \quad \text{for } j\neq p,p+1.\]
\end{lem}

\pf   Let $\lambda^{(1)}=(\lambda^{(1)}_1,\,\lambda^{(1)}_2,\ldots, \lambda^{(1)}_{N_1})$
be the first sub-overpartition of $\lambda$ in $\mathbb{F}^{(1)}_{N_1,\ldots,N_{k-1};i,p}$.  By definition, we see that
$\lambda_p^{(1)}$ is a non-overlined odd part, set $\lambda_p^{(1)}=2t+1$. Since  $\lambda \in \mathbb{F}^{(1)}_{N_1,\ldots,N_{k-1};i,p}$, we see that $2t+2$ does not occur in $\lambda$ when $|\lambda_{p-1}^{(1)}|\leq 2t-1$ and $\lambda_j^{(1)}$ is an  overlined odd part or a non-overlined even part, where $p+1\leq j\leq N_1$. When $1<p<N_1$, set $\lambda_{p+1}^{(1)}=\overline{2a+1}$ or $2a+2$,  it follows from the definition of \GG\ that  $a\geq t+1$.

For $1<p\leq N_1$, define $\mu=\Phi_{1,p}(\lambda)$, which can be obtained from $
\lambda$ by performing   the following two operations.

(1) Replace $\lambda_p^{(1)}=2t+1$ by   $\overline{2t+3}$.

(2) When $1<p<N_1$, replace $\lambda_{p+1}^{(1)}$ by    $2a+1$  (or  $\overline{2a+2}$) if $\lambda_{p+1}^{(1)}=\overline{2a+1}$ (or $2a+2$); When $p=N_1$, we shall do nothing.

Obviously, $|\mu|=|\lambda|+2$. We first prove that
  the parts from $\lambda$ in $\mu$ have the same marks as in $\lambda$ and the new generated parts $\overline{2t+3}$ and $2a+1$ (or $\overline{2a+2}$) in $\mu$ are   marked with 1.  This leads to    $\mu \in \mathbb{F}_{N_1,\ldots,N_{k-1};i}$.
  By the definition of $\mu$, it is  obvious that this assertion is true for the parts of size not exceed $2t+1$  in $\mu$. We next show that     the parts   $2t+2$ in $\mu$ have the same marks as in $\lambda$  and the new generated part $\overline{2t+3}$ in $\mu$ is marked with $1$.   There are two cases: if $\lambda^{(1)}_{p-1}=\overline{2t}$, or $2t$, or $\overline{2t+1}$, or $2t+1$, then  $\mu^{(1)}_{p-1}=\overline{2t}$, or $2t$, or $\overline{2t+1}$, or $2t+1$, it follows from the definition of \GG \ that  the parts  $2t+2$ in $\mu$ have the same marks as in $\lambda$  even if the part $\lambda^{(1)}_{p}=2t+1$ replaced by $\overline{2t+3}$ in $\mu$. Furthermore,  the new generated part $\overline{2t+3}$ in $\mu$ is marked with 1 since there is no $1$-marked  part of size  $2t+2$ in $\mu$. If $|\lambda^{(1)}_{p-1}|\leq 2t-1$, then   $2t+2$ does not occur in $\lambda$, and so neither in $\mu$.   It follows that the new generated part $\overline{2t+3}$ in $\mu$ is marked with $1$. Therefore, in either case,     the parts   $2t+2$ in $\mu$ have the same marks as in $\lambda$ and the new generated part $\overline{2t+3}$ in $\mu$ is marked with $1$.
    By the definition of  \GG, it is easy to see that the new generated part $2a+1$ (or $\overline{2a+2}$) in $\mu$ replacing  $\lambda^{(1)}_{p+1}$  is marked with 1, which has the same size with $\lambda^{(1)}_{p+1}$ in $\lambda$, and so   the  parts of size larger than $2t+2$ in $\mu$ have the same marks as in $\lambda$. Thus, we arrive at our assertion and prove that $\mu \in \mathbb{F}_{N_1,\ldots,N_{k-1};i}$.

Let  $\mu^{(1)}=(\mu^{(1)}_1,\,\mu^{(1)}_2,\ldots, \mu^{(1)}_{N_1})$ be the first sub-overpartition of $\mu$. It follows from the above proof   that
$\lambda_j^{(1)}=\mu_j^{(1)}$ for $j\neq p,p+1$. Furthermore, $\mu_{p}^{(1)}=\overline{2t+3}$, $\mu_{p+1}^{(1)}$ is a non-overlined odd part or  an overlined even part, and $\mu_j^{(1)}$ is an  overlined odd part or a non-overlined even part, where $p+2\leq j\leq N_1$.  Moreover,  if $|\mu_{p+1}^{(1)}|\geq 2t+5$,  then $2t+4$ does not occur in $\mu$. This is because that $2t+4$ does not occur in $\lambda$ when $|\lambda_{p+1}^{(1)}|\geq 2t+5$ and $\lambda_{p}^{(1)}=2t+1$.   Therefore,  $\mu \in \mathbb{\overline{F}}^{(1)}_{N_1,\ldots,N_{k-1};i,p}$. Furthermore, it is not difficult to show that  $\Phi_{1,p}$ is reversible. So, we conclude that $\Phi_{1,p}$ is a bijection between $\mathbb{F}^{(1)}_{N_1,\ldots,N_{k-1};i,p}$ and $\mathbb{\overline{F}}^{(1)}_{N_1,\ldots,N_{k-1};i,p}$. This completes the proof.  \qed

For example, for $p=3$, let
\begin{equation*}\label{example1}
GG(\lambda)=\begin{array}{ccc}\left[
\begin{array}{cccccccccccccccccccccc}
 &&&&&6\\
 &2&&4&&&\overline{7}\\
 \overline{1}&&{\bf3}^2&&&{\bf 6}
\end{array}
\right]&\begin{array}{c}3\\2\\1\end{array}\end{array}
\end{equation*}
be the \GG\ representation of $\lambda$ in $\mathbb{F}^{(1)}_{4,3,1;3,3}$. Applying the bijection $\Phi_{1,3}$ to $\lambda$, we get
\begin{equation*}\label{example1}
GG(\mu)=\begin{array}{ccc}\left[
\begin{array}{cccccccccccccccccccccc}
 &&&&&6\\
 &2&&4&&&\overline{7}\\
 \overline{1}&&3&&{\bf \overline{5}}&{\bf \overline{6}}
\end{array}
\right]&\begin{array}{c}3\\2\\1\end{array}\end{array},
\end{equation*}
which is in $\mathbb{\overline{F}}^{(1)}_{4,3,1;3,3}$.

\begin{lem}\label{theta2} For $1<p\leq N_1$,
there is a bijection $\Phi_{2,p}$ between $\mathbb{F}^{(2)}_{N_1,\ldots,N_{k-1};i,p}$ and  $\mathbb{\overline{F}}^{(2)}_{N_1,\ldots,N_{k-1};i,p}$. Furthermore, for $\lambda \in \mathbb{{F}}^{(2)}_{N_1,\ldots,N_{k-1};i,p}$ and $\mu=\Phi_{2,p}(\lambda)\in \overline{\mathbb{F}}^{(2)}_{N_1,\ldots,N_{k-1};i,p}$, we have
\[|\mu|=|\lambda|+2,\quad \text{and} \quad \lambda_j^{(1)}=\mu_j^{(1)} \quad \text{for } j\neq p,p+1.\]
\end{lem}

\pf Let $\lambda^{(1)}=(\lambda^{(1)}_1,\,\lambda^{(1)}_2,\ldots, \lambda^{(1)}_{N_1})$
be the first sub-overpartition of $\lambda$ in $\mathbb{F}^{(2)}_{N_1,\ldots,N_{k-1};i,p}$. By definition, we see that
$\lambda_p^{(1)}$ is a non-overlined odd part, set $\lambda_p^{(1)}=2t+1$. Note that $\lambda \in \mathbb{F}^{(2)}_{N_1,\ldots,N_{k-1};i,p}$, we see that  $|\lambda_{p-1}^{(1)}|\leq 2t-1$ and there is at least one part $2t+2$ in $\lambda$. Furthermore, $\lambda_{p+1}^{(1)}$ is an overlined odd part or a non-overlined even part, set $\lambda_{p+1}^{(1)}=\overline{2a+1}$ or $2a+2$. From the definition of \GG, we see that  $a\geq t+1$.

To define   $\Phi_{2,p}$, we  introduce an    index $r$ which is related to the marks of the parts $2t+2$ in $\lambda$. If there exists an integer $b$ such that there are  $b$-marked parts $2t$ and  $2t+2$ in $\lambda$, then set $r=b$; Otherwise, we define $r$ as the largest mark of the parts $2t+2$ in $\lambda$. Since $2t+2$ occurs in $\lambda$ and $\lambda_p^{(1)}=2t+1$, by the definition of  G\"ollnitz-Gordon marking,  we deduce that the marks of the parts $2t+2$ in $\lambda$ are larger than $1$. Hence  $r\geq2$.

For $1<p\leq N_1$, define $\mu=\Phi_{2,p}(\lambda)$, which  can be obtained from $\lambda$ by doing the following two operations.

(1) Replace $\lambda_p^{(1)}=2t+1$ by $2t+2$ and replace the $r$-marked part $2t+2$ in $\lambda$ by $\overline{2t+3}$.

(2) When $1<p<N_1$, replace $\lambda_{p+1}^{(1)}$ by    $2a+1$  (or  $\overline{2a+2}$) if $\lambda_{p+1}^{(1)}=\overline{2a+1}$ (or $2a+2$); When $p=N_1$, we shall do nothing.

 Obviously, $|\mu|=|\lambda|+2$. We first show that  $\mu$ is an overpartition in $ \mathbb{F}_{N_1,\ldots,N_{k-1};i}$. To this end,  we shall show that
  the parts from $\lambda$ in $\mu$ have the same marks as in $\lambda$ and the new generated parts ${2t+2}$ and $2a+1$ (or $\overline{2a+2}$) in $\mu$ are   marked with $1$, and the new generated part  $\overline{2t+3}$ in $\mu$  is   marked with $r$.   By the definition of $\mu$, it is  obvious that  the marks of parts of size not exceed $2t+1$ in $\mu$ are the same as  in $\lambda$. Note that $|\lambda_{p-1}^{(1)}|\leq 2t-1$,  $\lambda_p^{(1)}=2t+1$ and $|\lambda_{p+1}^{(1)}|>2t+1$, this implies that there are no  $1$-marked parts $\overline{2t}$, $2t$ and $\overline{2t+1}$ in $\lambda$, and there is only one part $2t+1$ in $\lambda$, that is $\lambda_p^{(1)}$. Therefore there are no  $1$-marked parts $\overline{2t}$, $2t$,  $\overline{2t+1}$  and   $2t+1$  in $\mu$. By the definition of  \GG, we see that the new generated part $2t+2$ in $\mu$ replacing  $\lambda_p^{(1)}=2t+1$  should be marked with 1 and  the parts   $2t+2$ from $\lambda$  in $\mu$ have  the same marks  as in $\lambda$.  We proceed to show that the new generated part $\overline{2t+3}$ in $\mu$  replacing the $r$-marked part $2t+2$ in $\lambda$ is also marked with $r$.
   By the definition of   \GG, we see that there are $f$-marked  $(2t+2)$'s in   $\lambda$ where $2\leq f\leq r-1$. It follows from the preceding proof  that there are also $f$-marked  $(2t+2)$'s in $\mu$, where $1\leq f\leq r-1$.  Therefore the new generated   part $\overline{2t+3}$ in $\mu$ replacing the $r$-marked $2t+2$ in $\lambda$ is also marked with $r$.
  Again, from the definition of  \GG, it is easy to see that the new generated part $2a+1$ (or $\overline{2a+2}$) in $\mu$ replacing $\lambda^{(1)}_{p+1}$   is marked with 1, which has the same size with $\lambda^{(1)}_{p+1}$ in $\lambda$, and so   the parts of size larger than $2t+1$ in $\mu$ have the same marks as in $\lambda$. Thus, we prove that $\mu$ is an overpartition in $ \mathbb{F}_{N_1,\ldots,N_{k-1};i}$.

Let  $\mu^{(1)}=(\mu^{(1)}_1,\,\mu^{(1)}_2,\ldots, \mu^{(1)}_{N_1})$ be the first sub-overpartition of $\mu$.  From the above proof,  we see that
$\lambda_j^{(1)}=\mu_j^{(1)}$ for $j\neq p,p+1$. Furthermore  $\mu_{p}^{(1)}=2t+2$, $\mu_{p+1}^{(1)}$ is a non-overlined odd part or an overlined even part, and $\mu_j^{(1)}$ is an  overlined odd part or a non-overlined even part, where $p+2\leq j\leq N_1$. This proves that $\mu$ is an overpartition in $ \mathbb{\overline{F}}_{N_1,\ldots,N_{k-1};{i,p}}$. Again by  the preceding proof, we also see that  there is an $r$-marked $\overline{2t+3}$ in $\mu$.  Furthermore, from the definition of $\lambda$, we see that    $2t+4$ does not occur in $\lambda$ if $|\lambda_{p+1}^{(1)}|> 2t+4$, otherwise, it contradicts to the assumptions that $\lambda_{p}^{(1)}=2t+1$ and $|\lambda_{p+1}^{(1)}|> 2t+4$. Hence,  we derive that $2t+4$ does not occur in $\mu$ if $|\mu_{p+1}^{(1)}|> 2t+4$. This proves that   $\mu$ is an overpartition in $\mathbb{\overline{F}}^{(2)}_{N_1,\ldots,N_{k-1};i,p}$. Moreover, it can be checked that
$\Phi_{2,p}$ is reversible. So we conclude that $\Phi_{2,p}$ is a bijection between $\mathbb{F}^{(2)}_{N_1,\ldots,N_{k-1};i,p}$ and  $\mathbb{\overline{F}}^{(2)}_{N_1,\ldots,N_{k-1};i,p}$.  This completes the proof. \qed

For example, for $p=3$, let
\begin{equation*}\label{example1}
GG(\lambda)=\begin{array}{ccc}\left[
\begin{array}{cccccccccccccccccccccc}
 &&&4&&&&{\bf 8}&&\\
 &2&&4&&&&8&&&\overline{11}\\
 \overline{1}&&&\overline{4}&&&{\bf 7}&&&{\bf 10}&&&\overline{13}
\end{array}
\right]&\begin{array}{c}3\\2\\1\end{array}\end{array}
\end{equation*}
be a \GG\ representation of $\lambda$ in $\mathbb{F}^{(2)}_{5,4,2;3,3}$. Applying the bijection $\Phi_{2,3}$ to $\lambda$, we get
\begin{equation*}\label{example1}
GG(\mu)=\begin{array}{ccc}\left[
\begin{array}{cccccccccccccccccccccc}
 &&&4&&&&&{\bf \overline{9}}&\\
 &2&&4&&&&8&&&\overline{11}\\
 \overline{1}&&&\overline{4}&&&&{\bf 8}&&{\bf \overline{10}}&&&\overline{13}
\end{array}
\right]&\begin{array}{c}3\\2\\1\end{array}\end{array},
\end{equation*}
which is in $\mathbb{\overline{F}}^{(2)}_{5,4,2;3,3}$.

\begin{lem}\label{theta3} For $1<p\leq N_1$,
there is a bijection $\Phi_{3,p}$ between $\mathbb{F}^{(3)}_{N_1,\ldots,N_{k-1};i,p}$ and  $\mathbb{\overline{F}}^{(3)}_{N_1,\ldots,N_{k-1};i,p}$. Furthermore, for $\lambda \in \mathbb{{F}}^{(3)}_{N_1,\ldots,N_{k-1};i,p}$ and $\mu=\Phi_{3,p}(\lambda)\in \overline{\mathbb{F}}^{(3)}_{N_1,\ldots,N_{k-1};i,p}$, we have
\[|\mu|=|\lambda|+2,\quad \text{and} \quad \lambda_j^{(1)}=\mu_j^{(1)} \quad \text{for } j\neq p,p+1.\]
\end{lem}

\pf Let $\lambda^{(1)}=(\lambda^{(1)}_1,\,\lambda^{(1)}_2,\ldots, \lambda^{(1)}_{N_1})$
be the first sub-overpartition of $\lambda$ in $\mathbb{F}^{(3)}_{N_1,\ldots,N_{k-1};i,p}$. By definition, we see that
$\lambda_p^{(1)}$ is an overlined even part, set  $\lambda_p^{(1)}=\overline{2t}$. Note that $\lambda \in\mathbb{F}^{(3)}_{N_1,\ldots,N_{k-1};i,p}$, so $\overline{2t+1}$ does not occur in $\lambda$ and   $\lambda_j^{(1)}$ is an  overlined odd part or a non-overlined even part, where  $p+1\leq j\leq N_1$. This implies that there are no parts of size $2t+1$ in $\lambda$. Set $\lambda_{p+1}^{(1)}=\overline{2a+1}$ or $2a+2$, where $a\geq t+1$.

To define $\Phi_{3,p}$, we introduce an index $r$ which is related to the marks of the parts  $2t$ in $\lambda$. There are three cases:  If $\lambda_{p-1}^{(1)}=\overline{2t-2}$, or $2t-2$, or $\overline{2t-1}$, or $2t-1$, then set $r=1$; If $|\lambda^{(1)}_{p-1}|\leq 2t-3$ and there exists an interger $b$ such that there are  $b$-marked  parts $2t-2$ and    $2t$ in $\lambda$, then set $r=b$; Otherwise,   set $r$ to be the largest mark of the parts of size  $2t$ in $\lambda$. Since $\lambda_p^{(1)}=\overline{2t}$,  we see that $r\geq 1$.

For $1<p\leq N_1$, define $\mu=\Phi_{3,p}(\lambda)$ which can be obtained from $\lambda$ by doing the following two operations:

(1) When $r=1$, replace $\lambda_p^{(1)}=\overline{2t}$  by    $2t+2$; When $r\geq 2$,  replace $\lambda_p^{(1)}=\overline{2t}$ by $2t$ and  replace the $r$-marked part $2t$ in $\lambda$ by   $2t+2$.

(2) When $1<p<N_1$, replace $\lambda_{p+1}^{(1)}$ by    $2a+1$  (or  $\overline{2a+2}$) if $\lambda_{p+1}^{(1)}=\overline{2a+1}$ (or $2a+2$); When $p=N_1$, we shall do nothing.

Obviously, $|\mu|=|\lambda|+2$. We first show that   $\mu$ is an overpartition in $\mathbb{F}_{N_1,\ldots,N_{k-1};i}$. We assert that the parts from $\lambda$ in $\mu$ have the same marks as in $\lambda$ and the new generated  part  $2a+1$ (or $\overline{2a+2}$) in $\mu$ is marked with $1$; When $r\geq 1$,   the new  generated part  ${2t+2}$ in $\mu$ replacing the $r$-marked part of size $2t$ is marked with  $r$; When $r\geq 2$, the new generated  part $2t$  replacing  $\lambda_p^{(1)}=\overline{2t}$   is marked with $1$.  From the construction of  $\mu$ and the definition of \GG, it is  obvious that  the marks of parts of size not exceed $2t-1$  in $\mu$ are the same as in $\lambda$. We proceed to show  that this assertion holds for   parts of size  not exceed $2t+2$ in $\mu$.    We consider the following two cases:

(1) When $r=1$, there are two subcases:

(1.1) When $\lambda_{p-1}^{(1)}=\overline{2t-2}$, or $2t-2$, or $\overline{2t-1}$, or $2t-1$, by the definition of   \GG, we see that   the parts $2t$  from $\lambda$ in $\mu$ have the same marks as   in $\lambda$  and  the new generated part $2t+2$ in $\mu$ replacing   $\lambda_p^{(1)}=\overline{2t}$ should be marked with 1. Since there are no parts of size $2t+1$ in $\lambda$,  from the construction of $\mu$, we see that there are no parts of size $2t+1$  in $\mu$. Note that the new generated part $2t+2$ in $\mu$ is marked with 1,  so we conclude  that the marks of   parts  $2t+2$ from $\lambda$ in $\mu$ are the same as in $\lambda$.

(1.2) When $|\lambda^{(1)}_{p-1}|\leq 2t-3$, and there is only one part of size   $2t$ in $\lambda$, that is, $\lambda_{p}^{(1)}=\overline{2t}$, from the construction of $\mu$, we see that there are no  parts of size $2t$ or $2t+1$ in $\mu$ and the new generated part $2t+2$ in $\mu$ replacing $\lambda_{p}^{(1)}=\overline{2t}$ is marked with $1$. Furthermore, the marks of the parts  $2t+2$ from $\lambda$ in $\mu$ stay the same as  in $\lambda$.

(2) When $r\geq2$, then either there are  $r$-marked parts  $2t-2$ and     $2t$ in $\lambda$, or $r$ is the largest mark of the parts   $2t$ in $\lambda$.  By the definition of \GG, we see that there are $f$-marked  parts of size $2t$ in $\lambda$, where $1\leq f\leq r$. It follows that the marks of the parts $2t+2$ in $\lambda$ are greater than $r$. Note that $|\lambda_{p-1}^{(1)}|\leq 2t-3$,  so the new generated part $2t$ in $\mu$ replacing  $\lambda_p^{(1)}=\overline{2t}$ is marked with $1$ and the marks of   parts $2t$ from $\lambda$ in $\mu$ are the same as  in $\lambda$. This means that there are $f$-marked    $2t$'s in $\mu$, where $1\leq f\leq r-1$. Since   there are no parts of size $2t+1$   in $\mu$,  we deduce that the new  generated part $2t+2$ in $\mu$ replacing the $r$-marked part $2t$ in $\lambda$ is marked with $r$ and the marks of  the parts $2t+2$ from $\lambda$ in $\mu$ are the same as in $\lambda$. Furthermore, $r$ is the smallest mark of the parts $2t+2$ in $\mu$.

  It remains to show that the assertion holds for  parts of size  greater than $2t+2$ in $\mu$.   By the definition of \GG, it is easy to see that the new generated $r$-marked part $2t+2$ in $\mu$ replacing the $r$-marked part of size $2t$ in $\lambda$ could affect the mark of the $r$-marked part of size $2t+3$ or $2t+4$ in $\mu$. Furthermore,  the new  generated part  $2a+1$ (or $\overline{2a+2}$) in $\mu$ replacing $\lambda^{(1)}_{p+1}$ is marked with $1$, which has the same size with $\lambda^{(1)}_{p+1}$ in $\lambda$. Hence it suffices to show that the $r$-marked part of size $2t+3$ or $2t+4$ from $\lambda$ in $\mu$ is also marked with $r$ in $\mu$.  We consider the following two cases:

(1) If $r=1$ and note that $\lambda_p^{(1)}=\overline{2t}$  and $|\lambda_{p+1}^{(1)}|\geq 2t+3$, then $\lambda^{(1)}_{p+1}=\overline{2t+3}$ (or  $2t+4$), this part    will become ${2t+3}$ (or $\overline{2t+4}$) in $\mu$, and by the definition of \GG, we see that ${2t+3}$ (or $\overline{2t+4}$) in $\mu$ is also marked with 1.

  (2) If $r\geq 2$, then we shall first show that there is no $r$-marked $\overline{2t+3}$ in $\mu$. By the definition of $\lambda$, we see that the marks of the parts of size $2t+2$ in $\lambda$ are greater than   $r$. It follows that $\overline{2t+3}$  in $\lambda$ is marked with 1, and so  there is no $r$-marked $\overline{2t+3}$ in $\mu$.  If there is an $r$-marked  $2t+4$ in $\lambda$, then there are $2$-marked,\ldots,$(r-1)$-marked parts $2t+4$ and a $1$-marked $\overline{2t+3}$ or $2t+4$ in $\lambda$.  Note that  $\lambda_{p+1}^{(1)}=\overline{2t+3}$ or $2t+4$ will become ${2t+3}$ or $\overline{2t+4}$ in $\mu$, which is also marked with 1 in $\mu$. Moreover, there is no $1$-marked $\overline{2t+2}$ or $2t+2$ or $\overline{2t+3}$ in $\mu$, and  $r$ is the smallest mark of the parts $2t+2$ in $\mu$, so the $r$-marked part $2t+4$ in $\lambda$ will also be marked with $r$ in $\mu$.

Hence we prove that the assertion holds for  parts of size  greater than $2t+2$ in $\mu$. So we conclude that $\mu$ is an overpartition in $\mathbb{F}_{N_1,\ldots,N_{k-1};i}$.

Let $\mu^{(1)}=(\mu^{(1)}_1,\,\mu^{(1)}_2,\ldots, \mu^{(1)}_{N_1})$ be the first sub-overpartition of $\mu$.  From the above proof, we see that
$\lambda_j^{(1)}=\mu_j^{(1)}$ for $j\neq p,p+1$. Furthermore,   $\mu_p^{(1)}$ is a non-overlined even part,  $\mu_{p+1}^{(1)}$ is a non-overlined odd part or an overlined even part, and $\mu_j^{(1)}$ is an  overlined odd part or a non-overlined even part, where $p+2\leq j\leq N_1$. Again, by  the above proof, we see that there is no overlined odd part of size $|\mu_p^{(1)}|+1$ in $\mu$. This proves that $\mu \in \mathbb{\overline{F}}^{(3)}_{N_1,\ldots,N_{k-1};i,p}$.

To prove that $\Phi_{3,p}$ is a bijection, we construct the inverse map $\Psi_{3,p}$ of
$\Phi_{3,p}$. Let  $\mu^{(1)}=(\mu^{(1)}_1,\,\mu^{(1)}_2,\ldots, \mu^{(1)}_{N_1})$ be the first sub-overpartition of $\mu$
 in  $\mathbb{\overline{F}}^{(3)}_{N_1,\ldots,N_{k-1};i,p}$. By definition, we see that $\mu_{p}^{(1)}$ is a non-overlined even part and $\mu_{p+1}^{(1)}$ is a non-overlined odd part or an overlined even part.  Let $\mu_{p}^{(1)}=2t$ and $\mu_{p+1}^{(1)}=2a+1$ or $\overline{2a+2}$. Note that $p>1$, so   $t\geq 2$. By the definition of \GG, we see that $a\geq t$. Since $\mu\in \mathbb{\overline{F}}^{(3)}_{N_1,\ldots,N_{k-1};i,p}$, $\overline{2t+1}$ does not occur in $\mu$.

To define  $\Psi_{3,p}$, we  define an index $r'$ related to the sizes of $\mu_{p}^{(1)}=2t$ and $\mu_{p+1}^{(1)}=2a+1$ or $\overline{2a+2}$.    If  $a=t$, or $a>t$ and   $2t+2$ does not occur in $\mu$, then set $r'=1$; If $a>t$ and     $2t+2$ occurs in $\mu$, then set $r'$ to be the smallest mark of the parts   $2t+2$ in $\mu$.

For $1<p\leq N_1$, define $\lambda=\Psi_{3,p}(\mu)$ which is obtained from $\mu$ by doing the following two operations.

(1) When $p=N_1$,  we shall do nothing; When $1<p<N_1$,    replace $\mu_{p+1}^{(1)}$ by $\overline{2a+1}$ (or $2a+2$)  if $\mu_{p+1}^{(1)}=2a+1$  (or $\overline{2a+2}$).

(2) When $r'=1$,  replace $\mu_{p}^{(1)}=2t$ by   $\overline{2t-2}$; When $r'\geq2$, replace the $r'$-marked  $2t+2$ in $\mu$ by $2t$ and replace $\mu_{p}^{(1)}=2t$ by  $\overline{2t}$ .

Obviously, $|\mu|=|\lambda|+2$.  It can be proved that $\lambda=\Psi_{3,p}(\mu)\in\mathbb{{F}}^{(3
  )}_{N_1,\ldots,N_{k-1};i,p}$ and $\Psi_{3,p}$ is the inverse map of $\Phi_{3,p}$. So, we conclude that $\Phi_{3,p}$ is a bijection between   $\mathbb{F}^{(3)}_{N_1,\ldots,N_{k-1};i,p}$ and  $\mathbb{\overline{F}}^{(3)}_{N_1,\ldots,N_{k-1};i,p}$. This completes the proof.  \qed

  For example, for $p=5$, let
\begin{equation*}\label{example1}
GG(\lambda)=\begin{array}{ccc}\left[
\begin{array}{cccccccccccccccccccccc}
 &&&4&&&&&&{\bf 10}&&&&14\\
 &2&&4&&&\overline{7}&&&10&&&&14\\
 \overline{1}&&3^2&&&6&&&&{\bf \overline{10}}&&&{\bf \overline{13}}
\end{array}
\right]&\begin{array}{c}3\\2\\1\end{array}\end{array}
\end{equation*}
be the \GG\ representation of $\lambda$ in $\mathbb{F}^{(3)}_{6,5,3;3,5}$. It can be checked that $r=3$. Applying the bijection $\Phi_{3,5}$ to $\lambda$, we get
\begin{equation*}\label{example1}
GG(\mu)=\begin{array}{ccc}\left[
\begin{array}{cccccccccccccccccccccc}
 &&&4&&&&&&&&{\bf 12}&&14\\
 &2&&4&&&\overline{7}&&&10&&&&14\\
 \overline{1}&&3^2&&&6&&&&{\bf 10}&&&{\bf 13}
\end{array}
\right]&\begin{array}{c}3\\2\\1\end{array}\end{array},
\end{equation*}
which is in $\mathbb{\overline{F}}^{(3)}_{6,5,3;3,5}$. Applying $\Psi_{3,5}$ to $\mu$, we see that $r'=3$ and $\Psi_{3,5}(\mu)=\lambda$.

For another example, for $p=5$, let
\begin{equation*}\label{example1}
GG(\lambda)=\begin{array}{ccc}\left[
\begin{array}{cccccccccccccccccccccc}
 &&&4&&&&&&&&&&14\\
 &2&&4&&&&8&&&&&&14\\
 \overline{1}&&&\overline{4}&&&7^2&{\bf \overline{8}}&&&&&{\bf \overline{13}}
\end{array}
\right]&\begin{array}{c}3\\2\\1\end{array}\end{array}
\end{equation*}
be the \GG\ representation of $\lambda$ in $\mathbb{F}^{(3)}_{6,4,2;3,5}$. Applying the bijection $\Phi_{3,5}$ to $\lambda$, we see that $r=1$ and
\begin{equation*}\label{example1}
GG(\mu)=\begin{array}{ccc}\left[
\begin{array}{cccccccccccccccccccccc}
 &&&4&&&&&&&&&&14\\
 &2&&4&&&&8&&&&&&14\\
 \overline{1}&&&\overline{4}&&&7^2&&&{\bf 10}&&&{\bf 13}
\end{array}
\right]&\begin{array}{c}3\\2\\1\end{array}\end{array},
\end{equation*}
which is in $\mathbb{\overline{F}}^{(3)}_{6,4,2;3,5}$. Applying $\Psi_{3,5}$ to $\mu$, we obtain that $r'=1$ and $\Psi_{3,5}(\mu)=\lambda$.

\begin{lem}\label{theta4} For $1<p\leq N_1$,
there is a bijection $\Phi_{4,p}$ between $\mathbb{F}^{(4)}_{N_1,\ldots,N_{k-1};i,p}$ and  $\mathbb{\overline{F}}^{(4)}_{N_1,\ldots,N_{k-1};i,p}$. Furthermore, for $\lambda \in \mathbb{{F}}^{(4)}_{N_1,\ldots,N_{k-1};i,p}$ and $\mu=\Phi_{4,p}(\lambda)\in \overline{\mathbb{F}}^{(4)}_{N_1,\ldots,N_{k-1};i,p}$, we have
\[|\mu|=|\lambda|+2,\quad \text{and} \quad \lambda_j^{(1)}=\mu_j^{(1)} \quad \text{for } j\neq p,p+1.\]
\end{lem}

\pf Let $\lambda^{(1)}=(\lambda^{(1)}_1,\,\lambda^{(1)}_2,\ldots, \lambda^{(1)}_{N_1})$
be the first sub-overpartition of $\lambda$ in $\mathbb{F}^{(4)}_{N_1,\ldots,N_{k-1};i,p}$. By definition, we see that
$\lambda_p^{(1)}$ is an overlined even part, set  $\lambda_p^{(1)}=\overline{2t}$.  So  $\overline{2t+1}$ occurs  in $\lambda$, assume that it is marked with  $s$ in $\lambda$. From the definition of \GG, it follows that $s\geq 2$. Note that  $\lambda_j^{(1)}$ is an  overlined odd part or a non-overlined even part, where $p+1\leq j\leq N_1$, set $\lambda_{p+1}^{(1)}=\overline{2a+1}$ or $2a+2$, by the definition of \GG, we see that $a\geq t+1$.

 To define $\Phi_{4,p}$, we also need to  use the index $r$ defined in the bijection $\Phi_{3,p}$. Recall that if $\lambda_{p-1}^{(1)}=\overline{2t-2}$, or $2t-2$, or $\overline{2t-1}$, or $2t-1$, then  $r=1$. If $|\lambda_{p-1}|\leq 2t-3$ and there exists an integer $b$ such that there are  $b$-marked parts $2t-2$ and  $2t$ in $\lambda$, then   $r=b$. Otherwise, set $r$ to be  the largest mark of the parts of size $2t$ in $\lambda$. By definition, we see that $s>r\geq1$.

For $1<p\leq N_1$, define $\mu=\Phi_{4,p}(\lambda)$ which can be obtained from $\lambda$ by doing the following two operations.

(1) When $r=1$, replace $\lambda_p^{(1)}=\overline{2t}$  by  $\overline{2t+1}$ and replace the $s$-marked $\overline{2t+1}$ in $\lambda$ by $2t+2$; When $r\geq 2$, first replace $\lambda_p^{(1)}=\overline{2t}$  by $2t$, and then replace the $r$-marked $2t$ in $\lambda$ by  $\overline{2t+1}$  and the   $s$-marked $\overline{2t+1}$  in $\lambda$ by $2t+2$.

(2) When $1<p<N_1$, replace $\lambda_{p+1}^{(1)}$ by    $2a+1$  (or  $\overline{2a+2}$) if $\lambda_{p+1}^{(1)}=\overline{2a+1}$ (or $2a+2$); When $p=N_1$, we shall do nothing.

Obviously, $|\mu|=|\lambda|+2$.  We first show that $\mu$ is an overpartition in $\mathbb{F}_{N_1,\ldots,N_{k-1};i}$. We   assert  that the parts from $\lambda$ in $\mu$ have the same marks as in $\lambda$ and the   new generated parts in $\mu$ replacing the parts in $\lambda$ have the same marks as   their original parts in $\lambda$.  By the definition of $\mu$, it is  obvious that  the marks of parts of size not exceed $2t-1$  in $\mu$ stay the same as   in $\lambda$. We proceed to show  that the marks of parts of size  $2t$ and $2t+1$ from $\lambda$ in $\mu$  are the same as in $\lambda$, the  new generated part $\overline{2t+1}$ replacing $\lambda^{(1)}_p=\overline{2t}$ is marked with $1$  when $r=1$,   and the  new generated part $2t$ replacing $\lambda^{(1)}_p=\overline{2t}$ is marked with $1$   and the  new generated part $\overline{2t+1}$ replacing the $r$-marked $2t$ in $\lambda$  is marked with $r$ when $r\geq 2$.  It should be mentioned that $2t+1$ does not occur in $\lambda$.   We consider the following two cases:

(1) When $r=1$, in this case, we see that $2t-2\leq |\lambda_{p-1}^{(1)}|\leq 2t-1$, or $|\lambda_{p-1}|\leq 2t-3$ and there is only one part of size $2t$ in $\lambda$. There are two subcases:

(1.1) When $\lambda_{p-1}^{(1)}=\overline{2t-2}$, or $2t-2$, or $\overline{2t-1}$, or $2t-1$,  by the definition of  \GG, we see that the marks of  parts $2t$ from $\lambda$ in $\mu$ are the same as in $\lambda$.  Note that $|\lambda_{p-1}^{(1)}|\leq 2t-1$, so the new generated   part $\overline{2t+1}$ in $\mu$ replacing $\lambda_p^{(1)}$ should be marked with 1.

(1.2) When $|\lambda_{p-1}|\leq 2t-3$ and there is only one part of size   $2t$ in $\lambda$, that is, $\lambda_{p}^{(1)}=\overline{2t}$,  there are no parts of size $2t$ in $\mu$, so the new generated  part $\overline{2t+1}$ in $\mu$ replacing $\lambda_{p}^{(1)}=\overline{2t}$ is marked with $1$.

(2) When $r\geq2$, either $|\lambda_{p-1}^{(1)}|\leq 2t-3$ and there are  $r$-marked parts $2t-2$ and   $2t$ in $\lambda$, or $r$ is the largest mark of the parts   $2t$ in $\lambda$. By the definition of \GG, we see that there are $f$-marked parts of size $2t$ in $\lambda$, where $1\leq f\leq r$.   Note that  $|\lambda_{p-1}|\leq 2t-3$, it follows that the new generated part $2t$ replacing $\lambda_p^{(1)}=\overline{2t}$ is marked with $1$ and the marks of  parts $2t$ from $\lambda$ in $\mu$ stay the same as in $\lambda$. Hence  there are   $f$-marked $2t$'s in $\mu$, where $1\leq f\leq r-1$. Therefore the new generated part  $\overline{2t+1}$ replacing the $r$-marked    $2t$ in $\lambda$ is marked with $r$ in $\mu$.

  Next, we show that the new generated   part $2t+2$  in $\mu$ replacing the $s$-marked $\overline{2t+1}$ in $\lambda$ is marked with $s$  and the marks of   parts $2t+2$ from $\lambda$ in $\mu$ are the same as in $\lambda$. Since there is an $s$-marked $\overline{2t+1}$ in $\lambda$, there are $f$-marked parts of size $2t$ in $\lambda$, where $1\leq f\leq s-1$. From the preceding proof and the  definition of $\mu$, it follows that there are $f$-marked   $2t$'s in $\mu$ where $1\leq f\leq s-1$ and $f\neq r$, and  there is an $r$-marked $\overline{2t+1}$ in $\mu$.  Hence the new generated part $2t+2$ in $\mu$ replacing the $s$-marked $\overline{2t+1}$ in $\lambda$  should be marked with $s$. Furthermore,   the marks of  parts $2t+2$ from $\lambda$ in $\mu$ are the same as in $\lambda$.

  It remains to show that the marks of  parts of size  greater than $2t+2$ in $\mu$ stay the same as in $\lambda$.   By the definition of  \GG, it is easy to see that the new generated $s$-marked $2t+2$ in $\mu$ replacing the $s$-marked $\overline{2t+1}$ in $\lambda$ only affect the mark of the $s$-marked part of size $2t+3$ or $2t+4$ in $\mu$.  Since $s\geq 2$, we see that there is no $s$-marked $\overline{2t+3}$ in $\lambda$,  and so there is no $s$-marked $\overline{2t+3}$ in $\mu$. Hence it suffices to show that the $s$-marked     $2t+4$ in $\lambda$ is also marked with $s$ in  $\mu$ even if the $s$-marked $\overline{2t+1}$ in $\lambda$ is replaced by the $s$-marked $2t+2$ in $\mu$.

 Note that there is an $s$-marked $\overline{2t+1}$ in $\lambda$, so there are a $1$-marked  $\overline{2t}$ and $f$-marked  $2t$'s in $\lambda$, where $2\leq f\leq s-1$. It follows that the marks of
 the parts $2t+2$ in $\lambda$ are greater than $s$. Hence we conclude that if there exists an $s$-marked  $2t+4$ in $\lambda$, then there are a $1$-marked $\overline{2t+3}$ or $2t+4$ and $f$-marked   $(2t+4)$'s in $\lambda$, where $2\leq f\leq s-1$.  Note that  $\lambda_{p+1}^{(1)}=\overline{2t+3}$ or $2t+4$ will become ${2t+3}$ or $\overline{2t+4}$ in $\mu$, which is also marked with 1 in $\mu$. Moreover, there is no $1$-marked $\overline{2t+2}$ or $2t+2$ or $\overline{2t+3}$ in $\mu$, and  $s$ is the least mark of the parts $2t+2$ in $\mu$, so by the definition of \GG, we see that the $s$-marked  $2t+4$ in $\lambda$ will also be marked with $s$ in   $\mu$.

Thus, we have shown that the marks of  parts  in $\mu$ are the same as the marks of their original parts in $\lambda$. Hence $\mu$ is an overpartition in $\mathbb{F}_{N_1,\ldots,N_{k-1};i}$.

Let $\mu^{(1)}=(\mu^{(1)}_1,\,\mu^{(1)}_2,\ldots, \mu^{(1)}_{N_1})$ be the first sub-overpartition of $\mu$.  It can be seen from  the above  proof that
$\lambda_j^{(1)}=\mu_j^{(1)}$ for $j\neq p,p+1$ and $\mu_{p}^{(1)}$ is an overlined odd part or a non-overlined even part. Furthermore, if $\mu_{p}^{(1)}$ is an overlined odd part, then  there is a non-overlined even part of size  $|\mu_p^{(1)}|+1$ in $\mu$, and $|\mu_{p+1}^{(1)}|\geq |\mu_p^{(1)}|+2$. If $\mu_p^{(1)}$ is a non-overlined even part, then  there are an overlined odd part of size $|\mu_p^{(1)}|+1$ and a non-overlined even part of size  $|\mu_p^{(1)}|+2$ in $\mu$, and  $|\mu_{p+1}^{(1)}|> |\mu_p^{(1)}|+2$.  Moreover, it is easy to see that $\mu_{p+1}^{(1)}$ is a non-overlined odd part or an overlined even part, and $\mu_j^{(1)}$ is an  overlined odd part or a non-overlined even part, where $p+2\leq j\leq N_1$.  This proves that $\mu$ is an overpartition in $ \mathbb{\overline{F}}^{(4)}_{N_1,\ldots,N_{k-1};i,p}$.

We proceed to construct the inverse map $\Psi_{4,p}$ of
$\Phi_{4,p}$, where  $1<p\leq N_1$. Let  $\mu^{(1)}=(\mu^{(1)}_1,\,\mu^{(1)}_2,\ldots, \mu^{(1)}_{N_1})$ be the first sub-overpartition of $\mu$
 in  $\mathbb{\overline{F}}^{(4)}_{N_1,\ldots,N_{k-1};i,p}$. By definition, we see that $\mu_{p}^{(1)}$ is an overlined odd part or a non-overlined even part and  $\mu_{p+1}^{(1)}$ is a non-overlined odd part or an overlined even part. If $\mu_{p}^{(1)}$ is an overlined odd part,     set $\mu_{p}^{(1)}=\overline{2t+1}$ and  $\mu_{p+1}^{(1)}= 2a+1 $ or $\overline{2a+2}$, then by definition, we see that $2t+2$ occurs in $\mu$ and $a\geq t+1$.  If $\mu_{p}^{(1)}$ is a non-overlined even part, set $\mu_{p}^{(1)}=2t$ and  $\mu_{p+1}^{(1)}= 2a+1 $ or $\overline{2a+2}$, then $\overline{2t+1}$ and $2t+2$ occur in $\mu$ and $a\geq t+1$.

Let  $r'$ be the mark of $\overline{2t+1}$ in $\mu$  and $s'$ be the smallest mark of the parts $2t+2$  in $\mu$. It follows from the definition of \GG \ that $s'>r'\geq 1$. Define $\lambda=\Psi_{4,p}(\mu)$ which is obtained from $\mu$ by doing the following two operations.

(1) When $p=N_1$,  we shall do nothing; When $p<N_1$,  replace $\mu_{p+1}^{(1)}$ by $\overline{2a+1}$ (or $2a+2$)    if $\mu_{p+1}^{(1)}=2a+1$  (or $\overline{2a+2}$).

(2) When $r'=1$,  replace $\mu_{p}^{(1)}=\overline{2t+1}$ by $\overline{2t}$ and replace the $s'$-marked $2t+2$ in $\mu$ by  $\overline{2t+1}$. When $r'\geq2$,  replace   $\mu_{p}^{(1)}={2t}$ by $\overline{2t}$, replace the $r'$-marked $\overline{2t+1}$ in $\mu$ by ${2t}$ and  replace the $s'$-marked $2t+2$ in $\mu$ by  $\overline{2t+1}$.

  It can be verified  that $\lambda=\Psi_{4,p}(\mu)\in\mathbb{{F}}^{(4)}_{N_1,\ldots,N_{k-1};i,p}$ and $|\mu|=|\lambda|+2$, and $\Psi_{4,p}$ is the inverse map of $\Phi_{4,p}$. So, we conclude that $\Phi_{4,p}$ is a bijection between $\mathbb{F}^{(4)}_{N_1,\ldots,N_{k-1};i,p}$ and  $\mathbb{\overline{F}}^{(4)}_{N_1,\ldots,N_{k-1};i,p}$. This completes the proof.   \qed

  For example, for $p=4$, let
\begin{equation*}\label{example1}
GG(\lambda)=\begin{array}{ccc}\left[
\begin{array}{cccccccccccccccccccccc}
 &&&&&&{\bf \overline{7}}&&&10\\
 &2&&&&6&&&&10\\
 \overline{1}&&&4&5&{\bf \overline{6}}&&&&{\bf 10}
\end{array}
\right]&\begin{array}{c}3\\2\\1\end{array}\end{array}
\end{equation*}
be the \GG\ representation of $\lambda$ in $\mathbb{F}^{(4)}_{5,3,2;3,4}$. Applying the bijection $\Phi_{4,4}$ to $\lambda$, we see that $r=1$, $s=3$, and
\begin{equation*}\label{example1}
GG(\mu)=\begin{array}{ccc}\left[
\begin{array}{cccccccccccccccccccccc}
 &&&&&&&{\bf 8}&&10\\
 &2&&&&6&&&&10\\
 \overline{1}&&&4&5&&{\bf \overline{7}}&&&{\bf \overline{10}}
\end{array}
\right]&\begin{array}{c}3\\2\\1\end{array}\end{array},
\end{equation*}
which is in $\mathbb{\overline{F}}^{(4)}_{5,3,2;3,4}$. Applying $\Psi_{4,4}$ to $\mu$, we have $r'=1$, $s'=3$, and  $\Psi_{4,4}(\mu)=\lambda$.

For example, for $p=3$, let
\begin{equation*}\label{example1}
GG(\lambda)=\begin{array}{ccc}\left[
\begin{array}{cccccccccccccccccccccc}
 &&&&&&{\bf \overline{7}}&&&&&&&14\\
 &2&&&&{\bf 6}&&&&10&&&\overline{13}\\
 \overline{1}&&3&&&{\bf \overline{6}}&&&{\bf \overline{9}}&&&12
\end{array}
\right]&\begin{array}{c}3\\2\\1\end{array}\end{array}
\end{equation*}
be the \GG\ representation of $\lambda$ in $\mathbb{F}^{(4)}_{5,4,2;3,3}$. Applying the bijection $\Phi_{4,3}$ to $\lambda$, we see that $r=2$, $s=3$, and
\begin{equation*}\label{example1}
GG(\mu)=\begin{array}{ccc}\left[
\begin{array}{cccccccccccccccccccccc}
 &&&&&&&{\bf 8}&&&&&&14\\
 &2&&&&&{\bf \overline{7}}&&&10&&&\overline{13}\\
 \overline{1}&&3&&&{\bf {6}}&&&{\bf {9}}&&&12
\end{array}
\right]&\begin{array}{c}3\\2\\1\end{array}\end{array},
\end{equation*}
which is in $\mathbb{\overline{F}}^{(4)}_{5,4,2;3,3}$. Applying $\Psi_{4,3}$ to $\mu$, we have $r'=2$, $s'=3$, and  $\Psi_{4,3}(\mu)=\lambda$.

We conclude this section by giving a proof of Lemma \ref{N-1-lem}. \vskip 0.2cm

{\bf Proof of Lemma \ref{N-1-lem}.} Supposed that $k\geq i\geq 1$, $N_1\geq N_2\geq\cdots \geq N_{k-1}\geq 0$ and $1< p\leq N_1$. From the definitions of  $\mathbb{F}^{(l)}_{N_1,\ldots,N_{k-1};i,p}$ and $\mathbb{\overline{F}}^{(l)}_{N_1,\ldots,N_{k-1};i,p}$, where $1\leq l\leq 4$,  we have
\[\mathbb{F}_{N_1,\ldots,N_{k-1};i,p}=\bigcup_{l=1}^{4}\mathbb{F}^{(l)}_{N_1,\ldots,N_{k-1};i,p}\]
and
\[\mathbb{\overline{F}}_{N_1,\ldots,N_{k-1};i,p}=\bigcup_{l=1}^{4}\mathbb{\overline{F}}^{(l)}_{N_1,\ldots,N_{k-1};i,p}.\]

Let $\lambda \in \mathbb{{F}}_{N_1,\ldots,N_{k-1};i,p}$, define
\[\mu=\Phi_p(\lambda)=\left\{\begin{array}{cc}\Phi_{1,p}(\lambda),&\text{ if }\lambda\in\mathbb{F}^{(1)}_{N_1,\ldots,N_{k-1};i,p};\\[5pt]
\Phi_{2,p}(\lambda),&\text{ if }\lambda\in\mathbb{F}^{(2)}_{N_1,\ldots,N_{k-1};i,p};\\[5pt]
\Phi_{3,p}(\lambda),&\text{ if }\lambda\in\mathbb{F}^{(3)}_{N_1,\ldots,N_{k-1};i,p};\\[5pt]
\Phi_{4,p}(\lambda),&\text{ if }\lambda\in\mathbb{F}^{(4)}_{N_1,\ldots,N_{k-1};i,p}.
\end{array}\right.\]
Combining Lemma \ref{theta1}, Lemma \ref{theta2}, Lemma \ref{theta3} and Lemma \ref{theta4}, we conclude that $\Phi_p$ is a bijection between  $\mathbb{F}_{N_1,\ldots,N_{k-1};i,p}$ and  $\mathbb{\overline{F}}_{N_1,\ldots,N_{k-1};i,p}$. Furthermore,    for $\lambda \in \mathbb{{F}}_{N_1,\ldots,N_{k-1};i,p}$ and $\mu=\Phi_{p}(\lambda)\in \overline{\mathbb{F}}_{N_1,\ldots,N_{k-1};i,p}$, we have
\[|\mu|=|\lambda|+2,\quad \text{and} \quad \lambda_j^{(1)}=\mu_j^{(1)}\quad \text{for} \quad j\neq p,p+1.\]

\qed

\section{Proof of Lemma \ref{N-2}}

 Let $\mathbb{D}_N$ denote the set of partitions $\eta=(\eta_1,\eta_2,\ldots, \eta_{\ell})$ with distinct negative odd parts which lay in $[1-2N,-1]$, that is,  $\eta_j$ is negative and odd for $1\leq  j\leq \ell$ and $1-2N\leq \eta_1<\eta_2<\cdots<\eta_\ell\leq -1$. It is clear to see that the generating function of partitions in  $\mathbb{D}_N$:
\[\sum_{\eta \in \mathbb{D}_{N}}q^{|\eta|}=(1+q^{1-2N})(1+q^{3-2N})\cdots (1+q^{-1})=(-q^{1-2N};q^2)_{N}.
\]

Hence Lemma \ref{N-2} is equivalent to the following combinatorial statement.

 \begin{thm}\label{N-2-COM} For $k\geq i\geq 1$ and $N_1\geq N_2\geq\cdots \geq N_{k-1}\geq 0$, there is  a bijection $\Theta$ between $\mathbb{{G}}_{N_1,\ldots,N_{k-1};i}$ and $\mathbb{D}_{N_1}\times\mathbb{{E}}_{N_1,\ldots,N_{k-1};i}$ such that   for $\mu \in \mathbb{{G}}_{N_1,\ldots,N_{k-1};i}$ and $\Theta(\mu)=(\eta,\nu)\in \mathbb{D}_{N_1}\times\mathbb{{E}}_{N_1,\ldots,N_{k-1};i}$, we have $|\mu|=|\eta|+|\nu|$.
\end{thm}

Let $\mu^{(r)}=(\mu^{(r)}_1,\mu^{(r)}_2,\ldots, \mu^{(r)}_{N_r})$ be the $r$-th sub-overpartition of  $\mu$ in  $\mathbb{G}_{N_1,\ldots,N_{k-1};i}$, where $1\leq r\leq k-1$. From the definition of $\mathbb{G}_{N_1,\ldots,N_{k-1};i}$, it is easy to see that $\mu_j^{(r)}$ is an overlined odd part or a non-overlined even part for $1\leq r\leq k-1$ and $1\leq j\leq N_r$.   Observe that $\mathbb{E}_{N_1,\ldots,N_{k-1};i}$ is the set of overpartitions in $\mathbb{G}_{N_1,\ldots,N_{k-1};i}$  for which there are no overlined odd parts, so the key point in the construction of the bijection $\Theta$  is   to  remove  all overlined odd parts  of an overpartition in $\mathbb{G}_{N_1,\ldots,N_{k-1};i}$ to get a new overpartition in $\mathbb{E}_{N_1,\ldots,N_{k-1};i}$.

Similarly to the bijection $\Phi$ in Section 4,  we will define three subsets $\mathbb{G}_{N_1,\ldots,N_{k-1};i,p}$, $\overline{\mathbb{G}}_{N_1,\ldots,N_{k-1};i,p}$ and $\overrightarrow{\mathbb{G}}_{N_1,\ldots,N_{k-1};i,p}$ of $\mathbb{G}_{N_1,\ldots,N_{k-1};i}$. Then we build a bijection $\Theta_p$ between $\mathbb{G}_{N_1,\ldots,N_{k-1};i,p}$ and  $\overline{\mathbb{G}}_{N_1,\ldots,N_{k-1};i,p}$ and a bijection $\Theta_{(p)}$ between $\mathbb{G}_{N_1,\ldots,N_{k-1};i,p}$ and  $\overrightarrow{\mathbb{G}}_{N_1,\ldots,N_{k-1};i,p}$. Similarly,   $\Theta_{(p)}$ can be obtained by successively using the bijection $\Theta_p$ and    plays a crucial role in the construction of the bijection $\Theta$ in Theorem \ref{N-2-COM}.

Let $\mu^{(1)}=(\mu^{(1)}_1,\mu^{(1)}_2,\ldots, \mu^{(1)}_{N_1})$ be  the first sub-overpartition of $\mu$ in $\mathbb{G}_{N_1,\ldots,N_{k-1};i}$.
To define the above three subsets of $\mathbb{G}_{N_1,\ldots,N_{k-1};i}$,   we  divide the parts in $\mu^{(1)}$ into two classes.  A part $\mu^{(1)}_j$ is   of type O if  $\mu^{(1)}_j$ is an overlined odd part or there is an overlined odd part of size $|\mu^{(1)}_j|+1$ in $\mu$, and   a part  $\mu^{(1)}_j$ is  of type E if  $\mu^{(1)}_j$ is a non-overlined even part and there is no overlined odd part of size $|\mu^{(1)}_j|+1$ in $\mu$. We say that two parts $\mu_{j_1}^{(1)}$ and $\mu_{j_2}^{(1)}$ are of  the same type if they are both of type O or they are both of type E.
For example, let
 \begin{equation*}
GG(\mu)=\left[\begin{array}{ccc}
\mu^{(3)}\\[3pt]
\mu^{(2)}\\[3pt]
\mu^{(1)}
\end{array}\right]=\begin{array}{ccc}\left[
\begin{array}{cccccccccccccccccccccc}
 &&&4&&&&&&10&&&&14&&\\
 &2&&&&6&&&\overline{9}&&&12&&&&\\
\overline{1}&&&4&&&&{8}&&&&12&&&&16
\end{array}
\right]&\begin{array}{c}3\\2\\1\end{array}\end{array}
\end{equation*}
be the \GG \ representation of $\mu$ in  $\mathbb{G}_{5,4,3;3}$.
 From definition, we see that the parts  $\mu_1^{(1)}=\overline{1}$ and $\mu_3^{(1)}=8$ are of type O, and the parts $\mu_2^{(1)}=4$, $\mu_4^{(1)}=12$ and $\mu_5^{(1)}=16$  are of  type E.

Let $k\geq i\geq 1$  and $N_1\geq N_2\geq\cdots \geq N_{k-1}\geq 0$ be given. For $1\leq p\leq N_1$, the subsets $\mathbb{G}_{N_1,\ldots,N_{k-1};i,p}$, $\overline{\mathbb{G}}_{N_1,\ldots,N_{k-1};i,p}$ and $\overrightarrow{\mathbb{G}}_{N_1,\ldots,N_{k-1};i,p}$ are described by using the first sub-overpartition $\mu^{(1)}=(\mu^{(1)}_1,\,\mu^{(1)}_2,\ldots, \mu^{(1)}_{N_1})$ of an overpartition $\mu$ in $\mathbb{G}_{N_1,\ldots,N_{k-1};i}$, where $\mu^{(1)}_1\leq \mu^{(1)}_2\leq \cdots\leq \mu^{(1)}_{N_1}$.

 $\lozenge$ Let $\mathbb{G}_{N_1,\ldots,N_{k-1};i,p}$ be the set of   overpartitions $\mu$ in $\mathbb{G}_{N_1,\ldots,N_{k-1};i}$ such that (1)
     $\mu_p^{(1)}$  is of type O for  $1\leq p\leq N_1$; (2) $\mu_j^{(1)}$  is of type  E, where $ p+1\leq j\leq N_1$.

$\lozenge$  Let $\mathbb{\overline{G}}_{N_1,\ldots,N_{k-1};i,p}$ be the set of   overpartitions $\mu$ in $\mathbb{G}_{N_1,\ldots,N_{k-1};i}$ such that (1)
     $\mu_p^{(1)}$  is of type E; (2) $\mu_{p+1}^{(1)}$  is of type  O; (3) $\mu_j^{(1)}$ is of type E, where  $p+2\leq j\leq N_1$.

$\lozenge$  Let $\overrightarrow{\mathbb{G}}_{N_1,\ldots,N_{k-1};i,p}$ be the set of   overpartitions $\mu$ in $\mathbb{G}_{N_1,\ldots,N_{k-1};i}$ such that $\mu_j^{(1)}$ is of type E, where $p\leq j\leq N_1$.

By definition, it is easy to see that for $1\leq p\leq N_1-2$
 \[\mathbb{\overline{G}}_{N_1,\ldots,N_{k-1};i,p}  \subseteq \mathbb{G}_{N_1,\ldots,N_{k-1};i,{p+1}} \subseteq  \overrightarrow{\mathbb{G}}_{N_1,\ldots,N_{k-1};i,p+2}.\]

 \begin{lem}\label{N-2-1-lem}
For $1\leq p \leq N_1$, there is a bijection $\Theta_p$ between $\mathbb{G}_{N_1,\ldots,N_{k-1};i,p}$ and   $\mathbb{\overline{G}}_{N_1,\ldots,N_{k-1};i,p}$.
Furthermore, for $\mu \in \mathbb{{G}}_{N_1,\ldots,N_{k-1};i,p}$ and $\nu=\Theta_p(\mu)\in
\mathbb{\overline{G}}_{N_1,\ldots,N_{k-1};i,p}$, we have
  $\mu_j^{(1)}$ and $\nu_j^{(1)}$ are of the same type for $ j\neq p,p+1$, and
$|\nu|=|\mu|+2-\delta_p^{N_1},$
where \[\delta_{p}^{N_1}=\begin{cases} 1, \text{ if } p=N_1;\\[3pt]
0, \text{ if } p\neq N_1.
\end{cases}
\]

 \end{lem}

 Applying in succession  the bijection in Lemma \ref{N-2-1-lem}  leads to a bijection between  $\mathbb{G}_{N_1,\ldots,N_{k-1};i,p}$ and $\overrightarrow{\mathbb{G}}_{N_1,\ldots,N_{k-1};i,p}$.

 \begin{lem}\label{N-2-2-lem}
For $1\leq p\leq N_1$, there is a bijection $\Theta_{(p)}$ between  $\mathbb{G}_{N_1,\ldots,N_{k-1};i,p}$ and  $\overrightarrow{\mathbb{G}}_{N_1,\ldots,N_{k-1};i,p}$. Furthermore, for  $\mu \in \mathbb{{G}}_{N_1,\ldots,N_{k-1};i,p}$ and $\nu=\Theta_{(p)}(\mu)\in\overrightarrow{\mathbb{G}}_{N_1,\ldots,N_{k-1};i,p}$, we have $\mu_j^{(1)}$ and $\nu_j^{(1)}$ are of the same type for $ j<p$ and
$|\nu|=|\mu|+2N_1-2p+1.$

 \end{lem}

  \pf Define $\Theta_{(p)}=\Theta_{N_1}\Theta_{N_1-1}\cdots \Theta_p$, by Lemma \ref{N-2-1-lem}, it is easy to verify that $\Theta_{(p)}$ is a bijection between   $\mathbb{G}_{N_1,\ldots,N_{k-1};i,p}$ and $\overrightarrow{\mathbb{G}}_{N_1,\ldots,N_{k-1};i,p}$ as desired. \qed

Before giving a proof of Lemma \ref{N-2-1-lem},  we first  give a proof of Lemma \ref{N-2} by successively applying the bijection $\Theta_{(p)}$ in Lemma \ref{N-2-2-lem}. \vskip 0.2cm

\noindent{\bf Proof of Lemma \ref{N-2}.} It is equivalent to prove Theorem \ref{N-2-COM}. Let
$\mu$ be  an overpartition  in $\mathbb{{G}}_{N_1,\ldots,N_{k-1};i}$.  We  aim to define $\Theta(\mu)=(\eta,\nu)$ such that $|\eta|+|\nu|=|\mu|$,  $\eta$ is a partition in $\mathbb{D}_{N_1}$ and $\nu$ is an overpartition in $\mathbb{{E}}_{N_1,\ldots,N_{k-1};i}$. We consider the following two cases.

Case 1. If there are no     overlined odd parts in $\mu$, then set $\nu=\mu$ and $\eta=\emptyset$.  It is easy to see that $\nu \in \mathbb{{E}}_{N_1,\ldots,N_{k-1};i}$ and $|\nu|=|\mu|$.

Case 2. If there are $s\geq 1$ overlined odd parts  in $\mu$, then there are    $s$ parts of type O   in the first sub-overpartition  of $\mu$. Note that if  there is an overlined odd part in $\mu$, say $\overline{2t+1}$, then it follows from the definition of \GG\ that there  exists a $1$-marked $\overline{2t+1}$ or  $2t$ in $\mu$.  So, each overlined  odd part in $\mu$ uniquely determines a part of type O in the first sub-overpartition  of $\mu$.

 Let  $\mu_{j_1}^{(1)}$, $\mu_{j_2}^{(1)}$, \ldots, $\mu_{j_{s-1}}^{(1)}$, $\mu_{j_s}^{(1)}$ be the parts of type O in the first sub-overpartition $\mu^{(1)}=(\mu_1^{(1)},\mu_2^{(1)},\ldots,\mu_{N_1}^{(1)})$ of $\mu$, where  $1\leq j_1<j_2< \cdots < j_s \leq N_1$. It is easy to see that $\mu \in \mathbb{{G}}_{N_1,\ldots,N_{k-1};i, j_{s}}$. Set
  \[\eta=(1-2(N_1-j_1+1), 1-2(N_1-j_2+1),\ldots,\ 1-2(N_1-j_s+1)).\]
   Obviously, $\eta \in {\mathbb{D}}_{N_1}$. The partition $\nu$ can be obtained from $\mu$ by employing the bijection in Lemma \ref{N-2-2-lem} $s$ times. We denote the intermediate partitions by $\gamma^{0},\ \gamma^{1},\ldots,\gamma^{s}$ with $\gamma^{0}=\mu$ and $\gamma^{s}=\nu$.
 For $1\leq b\leq s$,  the intermediate partition $\gamma^{b}$ can be obtained from  $\gamma^{b-1}$ by using $\Theta_{(j_{s-b+1})}$ in Lemma  \ref{N-2-2-lem}, that is,  for  $1\leq b\leq s$,
 \[\gamma^{b}=\Theta_{(j_{s-b+1})}(\gamma^{b-1}).\]
   Note that $\gamma^{0} \in \mathbb{{G}}_{N_1,\ldots,N_{k-1};i, j_s}$, so  by Lemma   \ref{N-2-2-lem}, we see that
   \[\gamma^{1} \in \mathbb{{G}}_{N_1,\ldots,N_{k-1};i, j_{s-1}}\text{ and }|\gamma^{1}|=|\mu|+2N_1-2j_s+1,\]
   and the first $(j_s-1)$ parts in the first sub-overpartitions  of $\gamma^{1}$ and $\gamma^{0}$ are of the same type.

  Successively employing Lemma \ref{N-2-2-lem}, we derive that for $1\leq b\leq s-1$,
   \[\gamma^{b} \in \mathbb{{G}}_{N_1,\ldots,N_{k-1};i, j_{s-p}}\text{ and }|\gamma^{b}|=|\mu|+\sum_{r=1}^b(2N_1-2j_{s-r+1}+1),\] and
   \[\gamma^{s} \in \overrightarrow{\mathbb{G}}_{N_1,\ldots,N_{k-1};i, j_{1}}\text{ and }|\gamma^{s}|=|\mu|+\sum_{r=1}^s(2N_1-2j_{s-r+1}+1).\]
   Furthermore, for $1\leq b\leq s$, the first $(j_{s-b+1}-1)$  parts in the first sub-overpartitions  of $\gamma^{b}$ and $\gamma^{0}$ are of the same type. From the definition of $\mu$,  the first $(j_1-1)$ parts in the first sub-overpartition of $\mu$ are of  type E, and by the definition of $\overrightarrow{\mathbb{G}}_{N_1,\ldots,N_{k-1};i, j_{1}}$,  we derive that   there are no  overlined odd parts in  $\gamma^{s}$. Hence
   \[\nu=\gamma^{s}\in \mathbb{E}_{N_1,\ldots,N_{k-1};i}\text{ and }|\nu|=|\mu|+\sum_{r=1}^s(2N_1-2j_{s-r+1}+1).\] It is easy to check that $|\eta|+|\nu|=|\mu|$.  Therefore $\Theta$ is well-defined.

 To prove that  $\Theta$ is a bijection, we shall give a brief description of the inverse map $\Lambda$ of $\Theta$.  Let $\nu$ be an overpartition in $\mathbb{E}_{N_1,\ldots,N_{k-1};i}$ and $\eta$ be a partition into distinct negative odd parts lying in $[1-2N_1,-1]$. We shall define  $\Lambda(\eta,\nu)=\mu$ such that $\mu$ is an overpartition in $\mathbb{G}_{N_1,\ldots,N_{k-1};i}$ and $|\eta|+|\nu|=|\mu|$.  There are two cases.

Case 1. If $\eta=\emptyset$, then set $\mu=\nu$. Note that $\mathbb{E}_{N_1,\ldots,N_{k-1};i}\subseteq \mathbb{G}_{N_1,\ldots,N_{k-1};i}$, so $\mu \in \mathbb{G}_{N_1,\ldots,N_{k-1};i}$ and there are no   overlined  odd parts in $\mu$.

Case 2. If $\eta\neq \emptyset$ and assume that
\[\eta=(1-2(N_1-j_1+1), 1-2(N_1-j_2+1),\ldots,\ 1-2(N_1-j_s+1)),\]
 where $1\leq j_1<j_2<\cdots<j_s\leq N_1$, then  $\mu$ can be recovered from $\nu$ by  using the bijection  in Lemma  \ref{N-2-2-lem} $s$ times. We denote the intermediate partitions by $\delta^{0},\,\delta^{1},\ldots,\delta^{s}$ with $\delta^{s}=\nu$ and $\delta^0=\mu$. For $1\leq b\leq s$, the intermediate partition $\delta^{b-1}$ can be obtained from  $\delta^{b}$ by using the   bijection $\Theta^{-1}_{(j_{s-b+1})}$  in Lemma  \ref{N-2-2-lem}, that is $\delta^{b-1}=\Theta^{-1}_{(j_{s-b+1})}(\delta^{b})$. Set $\mu=\delta^0$, by Lemma  \ref{N-2-2-lem},  we derive that $\mu$ is  an overpartition in $\mathbb{G}_{N_1,\ldots,N_{k-1};i}$ and $|\mu|=|\eta|+|\nu|$, and $\Lambda(\Theta(\mu))=\mu$ for any $\mu$ in $\mathbb{G}_{N_1,\ldots,N_{k-1};i}$. Hence  we conclude that $\Theta$ is a bijection between $\mathbb{{G}}_{N_1,\ldots,N_{k-1};i}$ and $\mathbb{D}_{N_1}\times\mathbb{{E}}_{N_1,\ldots,N_{k-1};i}$.
 This completes the proof of Theorem \ref{N-2-COM}, and hence Lemma \ref{N-2} is verified.   \qed

We proceed to give a proof of Lemma \ref{N-2-1-lem}. \vskip 0.2cm

 \noindent{\bf Proof of Lemma \ref{N-2-1-lem}.} Let $\mu^{(1)}=(\mu^{(1)}_1,\mu^{(1)}_2,\ldots, \mu^{(1)}_{N_1})$  be the first sub-overpartition of $\mu$ in  $\mathbb{G}_{N_1,\ldots,N_{k-1};i,p}$. By definition, we see that $\mu_p^{(1)}$ is of type O and $\mu_j^{(1)}$  is of type E, where $p+1\leq j\leq N_1$. If $\mu_p^{(1)}$ is an overlined odd part, then   set $\mu_p^{(1)}=\overline{2t+1}$; If $\mu_p^{(1)}$ is a non-overlined even part,   then  set $\mu_p^{(1)}=2t$, by definition of type O, we see that  there is an $s$-marked $\overline{2t+1}$ in $\mu$, where $s\geq 2$.

For $1\leq p\leq N_1$, define $\nu=\Theta_{p}(\mu)$ which can be obtained from $\mu$ as follows. There are three  cases.
\begin{itemize}
\item[Case 1] If $1\leq p<N_1$ and $\mu_p^{(1)}=\overline{2t+1}$, then $\mu_{p+1}^{(1)}$ is of type E, set $\mu_{p+1}^{(1)}=2b+2$, and it follows from the definition of \GG\ that $b\geq t+1$. There are two subcases.
     \begin{itemize}
   \item[Case 1.1]  If $b=t+1$, that is, $\mu_{p+1}^{(1)}=2t+4$,  then replace $\mu_p^{(1)}=\overline{2t+1}$ by  $2t+2$ and replace  $\mu_{p+1}^{(1)}=2t+4$ by $\overline{2t+5}$.
   \item[Case 1.2] If $b>t+1$, and set $r$ to be the largest mark of the parts $2b+2$ in $\mu$,   then replace $\mu_p^{(1)}=\overline{2t+1}$ by  $2t+2$ and replace the  $r$-marked $2b+2$ in $\mu$ by $\overline{2b+3}$.
     \end{itemize}

\item[Case 2] If   $1\leq p<N_1$ and $\mu_p^{(1)}=2t$, then $\mu_{p+1}^{(1)}$ is of type E, set $\mu_{p+1}^{(1)}=2b+2$ and $b\geq t+1$. Note that there is an $s$-marked $\overline{2t+1}$ in $\mu$. There are two subcases.
\begin{itemize}
   \item[Case 2.1]If there is  an  $s$-marked $2t+4$ in $\mu$, then
  replace the $s$-marked $\overline{2t+1}$ in $\mu$ by  $2t+2$ and replace the $s$-marked $2t+4$ in $\mu$  by $\overline{2t+5}$.

\item[Case 2.2] If  there is no    $s$-marked $2t+4$ in $\mu$, and set   $r$ to be the largest mark of the parts $2b+2$ in $\mu$, then
  replace the $s$-marked $\overline{2t+1}$ in $\mu$ by  $2t+2$ and replace the $r$-marked $2b+2$ in $\mu$ by  $\overline{2b+3}$.
  \end{itemize}

 \item[Case 3] If $p=N_1$, then there are two subcases.

\begin{itemize}
\item[Case 3.1] If $\mu_p^{(1)}=\overline{2t+1}$,  then replace $\mu_p^{(1)}=\overline{2t+1}$ by  $2t+2$.

\item[Case 3.2] If $\mu_p^{(1)}=2t$, and there is an $s$-marked $\overline{2t+1}$ in $\mu$, then
  replace the $s$-marked $\overline{2t+1}$ in $\mu$ by  $2t+2$.

\end{itemize}

\end{itemize}

Obviously, when $1\leq p<N_1$,  $|\nu|=|\mu|+2$, and when $p=N_1$, $|\nu|=|\mu|+1$. We next show that the parts from $\mu$ in $\nu$ have the same marks as in $\mu$ and the   new generated parts in $\nu$ replacing the parts in $\mu$ have the same marks as   their original parts in $\mu$. This implies that  $\nu\in\mathbb{{G}}_{N_1,\ldots,N_{k-1};i}$.  By the definition of $\nu$, it is  obvious that  the marks of parts of size not exceed $2t-1$  in $\nu$ stay the same as in $\mu$.  We now consider the marks of the new generated parts in $\nu$. There are two cases:

\begin{itemize}
\item If $\mu_p^{(1)}=\overline{2t+1}$, then  $|\mu_{p-1}^{(1)}|\leq 2t-1$, and so  there is no 1-marked $2t$ in $\nu$. Hence the new generated part $2t+2$ replacing $\mu_p^{(1)}=\overline{2t+1}$ in $\mu$ should be marked with $1$ in  $\nu$ and the parts $2t+2$ from $\mu$ in $\nu$ have the same marks as  in $\mu$. Thus, the marks of parts   $2t+2$  in $\nu$ stay the same as in $\mu$ when $\mu_{p}^{(1)}=\overline{2t+1}$ for $1\leq p\leq N_1$.

    For Case 1.1, since $\nu_{p}^{(1)}=2t+2$, it follows that the new generated part $\overline{2t+5}$ replacing $\mu_{p+1}^{(1)}=2t+4$ is   marked with $1$ in $\nu$.

    For Case 1.2, since $\nu_{p}^{(1)}=2t+2$ and $\mu_{p+1}^{(1)}=2b+2$, where $b>t+1$,  it follows from the definition of \GG \ that  the new generated part $\overline{2b+3}$ replacing the $r$-marked $2b+2$ in $\mu$ is   marked with $r$ in  $\nu$.

\item If $\mu_p^{(1)}=2t$,  then  there is an $s$-marked $\overline{2t+1}$  in $\mu$. It follows from the definition of \GG\  that
  there are  $f$-marked  $2t$'s   in $\mu$, where $1\leq f\leq s-1$, and hence the new generated part $2t+2$ replacing the $s$-marked $\overline{2t+1}$ in $\mu$ is also marked with $s$ in  $\nu$ and the marks of the parts $2t+2$ from $\mu$ in $\nu$ are the same as in $\mu$ when $\mu_{p}^{(1)}=2t$ for $1\leq p\leq N_1$.

 For Case 2.1,  note that there is an $s$-marked $2t+4$ in $\mu$, so there are $f$-marked  $(2t+4)$'s in $\mu$, where $1\leq f\leq s$, and  the parts $2t+4$ from $\mu$ in $\nu$ have the same marks as in $\mu$. Furthermore the new generated part $\overline{2t+5}$ replacing the $s$-marked $2t+4$ in $\mu$  is marked with $s$ in  $\nu$.

 For Case 2.2, note that there are $f$-marked $(2b+2)$'s in $\mu$, where $1\leq f\leq r$,  so we derive that the parts $2b+2$ from $\mu$ in $\nu$ have the same marks as in $\mu$ and the  new generated part $\overline{2b+3}$ in $\nu$ replacing the $r$-marked $2b+2$ in $\mu$ should be marked with  $r$ in $\nu$.

\end{itemize}

In any cases, we see that the marks of new generated parts in $\nu$ are the same as the marks of their original parts in $\mu$.   Furthermore,  the marks of the other parts from $\mu$ in $\nu$  are the same as in $\mu$. Hence   $\nu\in\mathbb{{G}}_{N_1,\ldots,N_{k-1};i}$.

  From the definition of $\nu$, it is easy to check that $\nu_{p}^{(1)}$ is of type E, $\nu_{p+1}^{(1)}$ is of type O, and $\nu_j^{(1)}$ is of type E, where $p+2\leq j\leq N_1$.   Hence, we deduce that $ \nu \in \mathbb{\overline{G}}_{N_1,\ldots,N_{k-1};i,p}$. Therefore,  $\Theta_p$ is well-defined. Furthermore,  $\mu_j^{(1)}$ and $\nu_j^{(1)}$ are of the same type for $ j\neq p,p+1$.

  To prove that $\Theta_{p}$ is a bijection, we give a brief description of the inverse map $\Lambda_{p}$ of
$\Theta_{p}$  for $1\leq p\leq N_1$. Let
 $\nu^{(1)}=(\nu^{(1)}_1,\nu^{(1)}_2,\ldots, \nu^{(1)}_{N_1})$  be the first sub-overpartition of $\nu$  in  $\mathbb{\overline{G}}_{N_1,\ldots,N_{k-1};i,p}$. By  definition,  $\nu_{p}^{(1)}$ is  of type E, $\nu_{p+1}^{(1)}$ is of type O, and $\nu_j^{(1)}$ is of type E, where $p+2\leq j\leq N_1$. For $1\leq p< N_1$, if $\nu_{p+1}^{(1)}$ is an overlined odd part, then set $\nu_{p+1}^{(1)}=\overline{2b+3}$; if $\nu_{p+1}^{(1)}$ is a non-overlined even part, then set $\nu_{p+1}^{(1)}=2b+2$, by the definition of  type O, we see that   there is an $r'$-marked $\overline{2b+3}$ in $\nu$, where $r'\geq 2$.

For $1\leq p\leq N_1$, define $\mu=\Lambda_{p}(\nu)$ which can be obtained from $\nu$ as follows. Here we set $\nu_p^{(1)}=2t+2$.  There are three cases.

 \begin{itemize}
 \item[Case 1] If $1\leq p<N_1$ and  $\nu_{p+1}^{(1)}=\overline{2b+3}$, then by the definition of \GG, we see that $t\leq b-1$. There are three subcases.
     \begin{itemize}
     \item[Case 1.1] If $t=b-1$, that is, $\nu_{p}^{(1)}=2b$, then replace $\nu_{p+1}^{(1)}=\overline{2b+3}$ by  $2b+2$ and replace  $\nu_{p}^{(1)}=2b$ by $\overline{2b-1}$.

    \item[Case 1.2] If $t<b-1$ and  $2t+4$ does not occur in $\nu$, then replace $\nu_p^{(1)}=2t+2$ by $\overline{2t+1}$ and replace $\nu_{p+1}^{(1)}=\overline{2b+3}$ by  $2b+2$.

     \item[Case 1.3] If  $t<b-1$ and $2t+4$  occurs  in $\nu$, set  $s'$ to be the smallest mark of the parts $2t+4$ in $\nu$, then  replace $\nu_{p+1}^{(1)}=\overline{2b+3}$ by  $2b+2$ and replace the $s'$-marked $2t+4$ in $\nu$ by $\overline{2t+3}$.
     \end{itemize}

 \item[Case 2] If $1\leq p<N_1$  and $\nu_{p+1}^{(1)}={2b+2}$, then it follows from  the definition of \GG\  that $t< b-1$ and there is an $r'$-marked  $\overline{2b+3}$  in   $\nu$. There are two subcases.

 \begin{itemize}
\item[Case 2.1] If  $2t+4$ does not occur   in $\nu$,  then replace the  $r'$-marked $\overline{2b+3}$ in $\nu$ by $2b+2$ and replace $\nu_p^{(1)}=2t+2$ by  $\overline{2t+1}$.

\item[Case 2.2] If  $2t+4$ occurs   in $\nu$ and set $s'$ to be the smallest mark of the parts $2t+4$ in $\nu$, then  replace the $r'$-marked $\overline{2b+3}$ in $\nu$ by $2b+2$ and
  replace the $s'$-marked $2t+4$ in $\nu$ by  $\overline{2t+3}$.

\end{itemize}

\item[Case 3] If $p=N_1$, then there are two subcases.

\begin{itemize}
\item[Case 3.1] If   $2t+4$ does not occur   in $\nu$,  then replace $\nu_p^{(1)}=2t+2$ by  $\overline{2t+1}$.

\item[Case 3.2] If   $2t+4$ occurs    in $\nu$, and set $s'$ to be the smallest mark of the parts $2t+4$ in $\nu$, then
  replace the $s'$-marked $2t+4$ in $\nu$ by  $\overline{2t+3}$.

\end{itemize}

 \end{itemize}

It can be verified  that $\Lambda_{p}(\nu)\in\mathbb{{G}}_{N_1,\ldots,N_{k-1};i,p}$   and $\Lambda_{p}$ is the inverse map of $\Theta_{p}$. So, we conclude that $\Theta_{p}$ is a bijection. \qed

We conclude this section by presenting the following three examples for the bijection $\Theta_{p}$ in  Lemma \ref{N-2-1-lem}.

(1) For $p=1$, let
\begin{equation*}\label{example1}
GG(\mu)=\begin{array}{ccc}\left[
\begin{array}{cccccccccccccccccccccc}
&&&&&&&8&&&&\\
&2&&&&&&8&&&&12\\
{\bf \overline{1}}&&&&&{\bf 6}&&&&10&&&&14
\end{array}
\right]&\begin{array}{c}3\\2\\1\end{array}\end{array}
\end{equation*}
be the \GG\ representation of $\mu$ in $\mathbb{{G}}_{4,3,1;3,1}$. Note that
$\mu^{(1)}_1=\overline{1}$ and $|\mu^{(1)}_2|=6>4$, which satisfy the conditions in Case 1.2 in the definition of $\Theta_p$, so $r=1$. Applying the bijection $\Theta_1$ to $\mu$, we get
\begin{equation*}\label{example1}
GG(\nu)=\begin{array}{ccc}\left[
\begin{array}{cccccccccccccccccccccc}
&&&&&&&8&&&&\\
&2&&&&&&8&&&&12\\
&{\bf 2}&&&&&{\bf \overline{7}}&&&10&&&&14
\end{array}
\right]&\begin{array}{c}3\\2\\1\end{array}\end{array},
\end{equation*}
which is in $\mathbb{\overline{G}}_{4,3,1;3,1}$. Note that
$\nu^{(1)}_1=2$ and $\nu^{(1)}_2=\overline{7}$ and   $4$ does not occur  in $\nu$, which satisfy the conditions in Case 1.2 in the definition of $\Lambda_p$. Applying $\Lambda_1$ to $\nu$, we recover $\mu$.

  (2) For $p=3$, let
\begin{equation*}\label{example1}
GG(\mu)=\begin{array}{ccc}\left[
\begin{array}{cccccccccccccccccccccc}
 &&&&&6&&&&&&12\\
 &2&&&\overline{5}&&&&{\bf \overline{9}}&&&{\bf 12}\\
 \overline{1}&&&4&&&&8&&&&12
\end{array}
\right]&\begin{array}{c}3\\2\\1\end{array}\end{array}
\end{equation*}
be the \GG\ representation of $\mu$ in $\mathbb{{G}}_{4,4,2;3,3}$. Note that
$\mu^{(1)}_3=8$, $\mu^{(2)}_3=\overline{9}$ and $\mu^{(2)}_4=12$, which satisfy the conditions in  Case 2.1 in the definition of $\Theta_p$, so $s=2$. Employing the bijection $\Theta_3$ to $\mu$, we get
\begin{equation*}\label{example1}
GG(\nu)=\begin{array}{ccc}\left[
\begin{array}{cccccccccccccccccccccc}
 &&&&&6&&&&&&12\\
 &2&&&\overline{5}&&&&&{\bf {10}}&&&{\bf \overline{13}}\\
 \overline{1}&&&4&&&&8&&&&12
\end{array}
\right]&\begin{array}{c}3\\2\\1\end{array}\end{array},
\end{equation*}
which is in $\mathbb{\overline{G}}_{4,4,2;3,3}$. Note that $\nu^{(1)}_3=8$, $\nu^{(2)}_3=10$, $\nu^{(1)}_4={12}$ and $\nu^{(2)}_4=\overline{13}$, which satisfy the conditions in Case 2.2 in the definition of $\Lambda_p$, so $r'=2$.  Applying $\Lambda_3$ to $\nu$, we recover $\mu$.

 (3) For $p=4$, let
\begin{equation*}\label{example1}
GG(\mu)=\begin{array}{ccc}\left[
\begin{array}{cccccccccccccccccccccc}
 &&&&&6&&&&&&12\\
 &2&&&&6&&&&&{\bf \overline{11}}&\\
 \overline{1}&&&4&&&\overline{7}&&&10
\end{array}
\right]&\begin{array}{c}3\\2\\1\end{array}\end{array}
\end{equation*}
be the \GG\ representation of $\mu$ in $\mathbb{G}_{4,3,2;3,4}$. Note that $\mu_4^{(1)}=10$ and $\mu_3^{(2)}=\overline{11}$, which satisfy the conditions in Case 3.2 in the definition of $\Theta_p$, so $s=2$.  Applying the bijection $\Theta_{4}$ to $\mu$, we get
\begin{equation*}\label{example1}
GG(\nu)=\begin{array}{ccc}\left[
\begin{array}{cccccccccccccccccccccc}
 &&&&&6&&&&&&12\\
 &2&&&&6&&&&&&{\bf 12}\\
 \overline{1}&&&4&&&\overline{7}&&&10
\end{array}
\right]&\begin{array}{c}3\\2\\1\end{array}\end{array},
\end{equation*}
which is in $\mathbb{\overline{G}}_{4,3,2;3,4}$. Note that $\nu_4^{(1)}=10$ and $\nu^{(2)}_3=12$, which satisfy the conditions in Case 3.2 in the definition $\Lambda_p$, so $s'=2$.  Applying $\Lambda_{4}$ to $\nu$, we recover $\mu$.

\section{Proof of Theorem \ref{Gollnitz-odd-1}}

 In this section, we complete the proof of Theorem \ref{Gollnitz-odd-1}. Using Lemma \ref{N-1} and Lemma \ref{N-2}, we   first give a proof of the formula for the generating function of ${F}_{k,i}(m,n)$   in Theorem \ref{thm-2N-1}.

 \vskip 0.2cm

\noindent{\bf Proof of Theorem \ref{thm-2N-1}.} First, we derive the following formula for the generating function  of the number of overpartitions $\lambda$ in $\mathbb{E}_{N_1,\ldots,N_{k-1};i}$ with the aid of  the identity \eqref{GENREATING-EEK} due to Kur\c{s}ung\"{o}z.
\begin{equation}\label{GENREATING-EE}
\sum_{\lambda\in\mathbb{E}_{N_1,\ldots,N_{k-1};i}}
q^{|\lambda|}=
\frac{
q^{2(N^2_1+\cdots+N^2_{k-1}+N_{i}+\cdots+N_{k-1})}
}
{(q^2;q^2)_{N_1-N_2}\cdots(q^2;q^2)_{N_{k-2}-N_{k-1}}
(q^{2};q^{2})_{N_{k-1}}},
\end{equation}

 Recall that  $\mathbb{B}_{N_1,\ldots,N_{k-1};i}$ is the set of ordinary partitions $\eta$ for which
\begin{equation}\label{pf-thm3.4-2}
f_1(\eta)\leq i-1\quad \text{and} \quad f_{t}(\eta)+f_{t+1}(\eta)\leq k-1
\end{equation}
such that there are $N_r$ $r$-marked parts in the  Gordon  marking of $\eta$ for $1\leq r\leq k-1$. To prove Theorem \ref{R-R-G-left-e} combinatorially, Kur\c{s}ung\"{o}z \cite{Kursungoz-2010} established the following formula for the generating function of the number of ordinary partitions in $\mathbb{B}_{N_1,\ldots,N_{k-1};i}$.

\begin{equation}\label{GENREATING-EEK-6}
\sum_{\eta\in\mathbb{B}_{N_1,\ldots,N_{k-1};i}}
q^{|\eta|}=
\frac{
q^{N^2_1+\cdots+N^2_{k-1}+N_{i}+\cdots+N_{k-1}}
}
{(q;q)_{N_1-N_2}\cdots(q;q)_{N_{k-2}-N_{k-1}}
(q;q)_{N_{k-1}}}.
\end{equation}

From the definitions in Section 3, we see that  $\mathbb{E}_{N_1,\ldots,N_{k-1};i}$ is also the set of ordinary partitions $\lambda$ without odd parts for which
\begin{equation}\label{pf-thm3.4-1}
f_2(\lambda)\leq i-1\quad \text{and} \quad f_{2t}(\lambda)+f_{2t+2}(\lambda)\leq k-1
\end{equation}
such that there are $N_r$ $r$-marked parts in the  G\"ollnitz-Gordon marking of $\lambda$ for $1\leq r\leq k-1$.

To show \eqref{GENREATING-EE}, we aim to build a bijection $\phi$ between  $\mathbb{E}_{N_1,\ldots,N_{k-1};i}$   and   $\mathbb{B}_{N_1,\ldots,N_{k-1};i}$ such that for $\lambda \in \mathbb{E}_{N_1,\ldots,N_{k-1};i}$ and $\eta=\phi(\lambda)\in \mathbb{B}_{N_1,\ldots,N_{k-1};i}$, we have $|\lambda|=2|\eta|$.

In terms of generating functions,  we have
\begin{eqnarray}\label{GENREATING-EEK-aaa}
& &\sum_{\lambda\in\mathbb{E}_{N_1,\ldots,N_{k-1};i}}
q^{|\lambda|}= \sum_{\eta\in\mathbb{B}_{N_1,\ldots,N_{k-1};i}}
q^{2|\eta|}.
\end{eqnarray}
Let $\lambda=(\lambda_1,\lambda_2,\ldots,\lambda_\ell)$  be a partition in $\mathbb{E}_{N_1,\ldots,N_{k-1};i}$, where $\lambda_1\leq \lambda_2 \leq \cdots \leq \lambda_\ell$,  we see that $\lambda_j$ is a non-overlined even part for $1\leq j\leq \ell$ and $\lambda$ satisfies  \eqref{pf-thm3.4-1}. Define
\[\eta=\phi(\lambda)=\left(\frac{\lambda_1}{2},\frac{\lambda_2}{2},\ldots,\frac{\lambda_\ell}{2}\right).\]
Clearly, $|\lambda|=2|\eta|$. Furthermore, $f_{t}(\eta)=f_{2t}(\lambda)$ for $t\geq 1$, which implies that $\eta$ satisfies \eqref{pf-thm3.4-2}. Hence it remains to show that the there are $N_r$ $r$-marked parts in the Gordon marking of $\eta$ for $1\leq r\leq k-1$.

By the definition of \GG, we see that the \GG \ of $\lambda=(\lambda_1,\lambda_2,\ldots,\lambda_\ell)$, where  $\lambda_j$ is a non-overlined even part for $1\leq j\leq \ell$   can be described as follows:  First, $\lambda_1$ is marked with 1,  and for $p>1$,  assume that the part $\lambda_j$ for $j<p$  has been assigned a mark. Then $\lambda_p$ is marked with   the least positive integer that is not used to mark the parts   $\lambda_j$ with $\lambda_p-\lambda_j\leq 2$  for $j<p$.  For example, the  G\"ollnitz-Gordon marking  of $\lambda=(2,2,4,4,4,6,8,10,10,12,12,12)$ is
\[GG(\lambda)=(2_1,2_2,4_3,4_4,4_5,6_1,8_2,10_1,10_3,12_2,12_4,12_5).\]
We now consider the Gordon marking of $\eta=(\eta_1,\eta_2,\ldots, \eta_\ell)$ where $\eta_1\leq \eta_2 \leq \cdots \leq \eta_\ell$. By definition, we see that $\eta_1$ is marked with 1,  and for $p>1$,  assume that the part $\eta_j$   has been assigned  a mark for $j<p$. Then $\eta_p$ is marked with   the least positive integer that is not used to mark the parts   $\eta_j$ with $\eta_p-\eta_j\leq 1$  for $j<p$. Since $\eta_j=\lambda_j/2$ for $1\leq j\leq \ell$,   it  can be checked that  the mark of $\eta_j$ in the Gordon marking of $\eta$  is the same as the mark of $\lambda_j$ in the G\"ollnitz-Gordon marking of $\lambda$ for $1\leq j\leq \ell$. For example, the  Gordon marking  of $\eta=\phi(\lambda)=(1,1,2,2,2,3,4,5,5,6,6,6)$ is
\[G(\eta)=(1_1,1_2,2_3,2_4,2_5,3_1,4_2,5_1,5_3,6_2,6_4,6_5).\]
Hence there are $N_r$ $r$-marked parts in the Gordon marking of $\eta$ for $1\leq r\leq k-1$, and so $\eta \in\mathbb{B}_{N_1,\ldots,N_{k-1};i}$. Furthermore, it is easy to see that this process is reversible. Therefore,  we conclude that $\phi$ is a bijection between  $\mathbb{E}_{N_1,\ldots,N_{k-1};i}$   and   $\mathbb{B}_{N_1,\ldots,N_{k-1};i}$, and   \eqref{GENREATING-EEK-aaa} holds.  Substituting \eqref{GENREATING-EEK-6} into \eqref{GENREATING-EEK-aaa}, we obtain \eqref{GENREATING-EE}.

Substituting \eqref{GENREATING-EE} into the relation \eqref{N2e} in Lemma \ref{N-2}, we   obtain the following generating function of the number of   overpartitions in $\mathbb{G}_{N_1,\ldots,N_{k-1};i}$.
\begin{equation}\label{N2ee}
\sum_{\mu\in\mathbb{G}_{N_1,\ldots,N_{k-1};i}}
q^{|\mu|}=\frac{(-q^{1-2N_1};q^2)_{N_1}
q^{2(N^2_1+\cdots+N^2_{k-1}+N_{i}+\cdots+N_{k-1})}
}
{(q^2;q^2)_{N_1-N_2}\cdots(q^2;q^2)_{N_{k-2}-N_{k-1}}
(q^{2};q^{2})_{N_{k-1}}}.
\end{equation}
Plugging \eqref{N2ee} into the relation \eqref{N1} in Lemma \ref{N-1}, we obtain \eqref{N1e}, which  yields  the generating function of ${F}_{k,i}(m,n)$  in Theorem \ref{thm-2N-1}. Thus we complete the proof. \qed

By Theorem \ref{thm-2N-1} and Lemma \ref{2N1}, we obtain the following generating function of $H_{k,i}(m,n)$.

\begin{thm}\label{GKI-1} For $k\geq i\geq1,$
\begin{equation}\label{GKI1}
\begin{split}
& \displaystyle\sum_{m,n\geq0}H_{k,i}(m,n)x^mq^n\\
&=\sum_{N_{1}\geq \cdots\geq N_{k-1}\geq 0}\frac{(-q^{2-2N_1};q^2)_{N_1-1}(-q^{1-2N_1};q^2)_{N_1}
q^{2(N^2_1+\cdots+N^2_{k-1}+N_{i+1}+\cdots+N_{k-1})}
x^{N_1+\cdots+N_{k-1}}}
{(q^2;q^2)_{N_1-N_2}\cdots(q^2;q^2)_{N_{k-2}-N_{k-1}}
(q^{2};q^{2})_{N_{k-1}}}.
\end{split}
\end{equation}
\end{thm}

\pf From the relation \eqref{2N-3} in Lemma \ref{2N1}, we deduce that for $1\leq i\leq k-1,$
\begin{equation}\label{GKI2}
\begin{split}
& \displaystyle\sum_{m,n\geq0}H_{k,i}(m,n)x^mq^n\\
&=\displaystyle\sum_{m,n\geq0}F_{k,i+1}(m,n)x^mq^n\\
&=\sum_{N_{1}\geq \cdots\geq N_{k-1}\geq 0}\frac{(-q^{2-2N_1};q^2)_{N_1-1}(-q^{1-2N_1};q^2)_{N_1}
q^{2(N^2_1+\cdots+N^2_{k-1}+N_{i+1}+\cdots+N_{k-1})}
x^{N_1+\cdots+N_{k-1}}}
{(q^2;q^2)_{N_1-N_2}\cdots(q^2;q^2)_{N_{k-2}-N_{k-1}}
(q^{2};q^{2})_{N_{k-1}}}.
\end{split}
\end{equation}
For $i=k$, from \eqref{2N-4} in Lemma \ref{2N1},  it follows that
\[\sum_{m,n\geq0}H_{k,k}(m,n)x^mq^n
=\sum_{m,n\geq0}F_{k,1}(m,n)(xq^{-2})^mq^n.\]
Using the generating function of $F_{k,1}(m,n)$, we obtain
\begin{equation}\label{GKI3}
\begin{split}
& \displaystyle\sum_{m,n\geq0}H_{k,k}(m,n)x^mq^n\\
&=\sum_{N_{1}\geq \cdots\geq N_{k-1}\geq 0}\frac{(-q^{2-2N_1};q^2)_{N_1-1}(-q^{1-2N_1};q^2)_{N_1}
q^{2(N^2_1+\cdots+N^2_{k-1})}
x^{N_1+\cdots+N_{k-1}}}
{(q^2;q^2)_{N_1-N_2}\cdots(q^2;q^2)_{N_{k-2}-N_{k-1}}
(q^{2};q^{2})_{N_{k-1}}}.
\end{split}
\end{equation}
Observe that the above formula \eqref{GKI2} for $1\leq i\leq k-1$ and \eqref{GKI3} for $i=k$ take the same form as   in Theorem \ref{GKI-1}. Thus, we complete the proof of  Theorem \ref{GKI-1}.  \qed

We conclude this paper with the proof of Theorem \ref{Gollnitz-odd-1}.

\noindent{\bf Proof of Theorem \ref{Gollnitz-odd-1}.} Substituting \eqref{thm-2N1} and \eqref{GKI1} into the relation \eqref{2N-2}, we obtain
\begin{equation*}
\begin{split}
& \displaystyle\sum_{m,n\geq0}O_{k,i}(m,n)x^mq^n\\
&=\displaystyle\sum_{m,n\geq0}F_{k,i}(m,n)x^mq^n+\displaystyle\sum_{m,n\geq0}H_{k,i}(m,n)x^mq^n\\
&=\sum_{N_{1}\geq \cdots\geq N_{k-1}\geq 0}\frac{(-q^{2-2N_1};q^2)_{N_1-1}(-q^{1-2N_1};q^2)_{N_1}
q^{2(N^2_1+\cdots+N^2_{k-1}+N_{i+1}+\cdots+N_{k-1})}(1+q^{2N_i})
x^{N_1+\cdots+N_{k-1}}}
{(q^2;q^2)_{N_1-N_2}\cdots(q^2;q^2)_{N_{k-2}-N_{k-1}}
(q^{2};q^{2})_{N_{k-1}}},
\end{split}
\end{equation*}
which is \eqref{Gollnitz-odd1}. This completes the proof.   \qed

 \vskip 0.2cm
\noindent{\bf Acknowledgments.} This work
was supported by the 973 Project, the PCSIRT Project and the National Science Foundation of China. We are greatly indebted to referees for their helpful suggestions that improved the presentation of this paper.

\end{document}